\documentclass[10pt,reqno,a4paper]{amsart}
\usepackage{euscript}
\usepackage{amssymb}
\usepackage{breqn}
\usepackage{graphicx}
\usepackage{color}
\usepackage{overpic}
\usepackage{pdflscape}
\renewcommand{\epsilon}{\varepsilon}
\usepackage{float}
\usepackage{soul}
\usepackage{hyperref}
\usepackage{epstopdf} 
\usepackage{caption}
\usepackage{amsthm,amssymb,amsmath,amsfonts,amscd}
\usepackage{graphicx,epsfig,latexsym}
\usepackage{cite,calc,xcolor,subfigmat,mathrsfs,dsfont,epstopdf}
\usepackage{tabularx}
\usepackage{circuitikz}

\newcommand{\RR}{\mathbb{R}}
\newcommand{\NN}{\mathbb{N}}

\newcommand{\0}{\mathbf{0}}

\newtheorem{thm}{Theorem}[section]

\newtheorem{lem}[thm]{Lemma}
\newtheorem{prop}[thm]{Proposition}
\theoremstyle{definition}
\newtheorem{defn}[thm]{Definition}

\newtheorem{rem}[thm]{Remark}

\numberwithin{equation}{section}
\def\Blem {\begin{lem}}
\def\Elem {\end{lem}}
\def\be {\begin{equation}}
\def\ee {\end{equation}}
\def\ba {\begin{eqnarray}}
\def\ea {\end{eqnarray}}
\def\bes {\begin{equation}}
\def\ees {\end{equation}}
\def\bas {\begin{eqnarray}}
\def\eas {\end{eqnarray}}
\def\bpr {\begin{proof}}
\def\epr {\end{proof}}

\begin{document}
 \title[FitzHugh-Nagumo system] {FitzHugh-Nagumo equation:  bifurcations, \\ slow-fast system and  dynamics near infinity\\ }
 
\maketitle

\begin{center}
\author{Alexandre A. P. Rodrigues$^{1}$ 
and Nasrin Sadri$^{2,\dag}$}\footnote{$\dag$ Corresponding author}\\\vspace*{0.5cm}

{$^1$ Lisbon School of Economics and Management, Lisbon University, \\
Centro de Matem\'atica Aplicada \`a Previs\~ao e Decis\~ao Econ\'omica, \\
Rua do Quelhas 6, 1200-781 Lisboa, Portugal}, 
\email{arodrigues@iseg.ulisboa.pt}

{$^2$ School of Mathematics, Statistics and Computer Science, College of Science, University of Tehran, Tehran, Iran and School of Mathematics, Institute for Research in Fundamental Sciences (IPM), P. O. Box 19395-5746, Tehran, Iran} \\ \email{n.sadri@ut.ac.ir}
\end{center}

\allowdisplaybreaks

\begin{abstract}
In this paper, we present a qualitative and bifurcation analysis of the three-parameter FitzHugh-Nagumo system and its compactified formulation. The study is structured according to three parameter-dependent regimes, for which the associated phase portraits are characterized.
In one of these regimes, the system exhibits a double-zero bifurcation with $\mathbb{Z}_2$--symmetry
 corresponding to a codimension-two degeneracy. We compute explicit bifurcation and transition curves arising in the unfolding of this singularity, including Pitchfork, Hopf, Belyakov, and double homoclinic bifurcations, and we construct the corresponding bifurcation diagrams. 
The local bifurcation analysis is linked to the slow-fast structure of the model, highlighting the emergence of canard solutions and their role in organizing the dynamics.

In addition, we analyse the global behavior of the system by studying the dynamics at infinity through compactification techniques. The results obtained complement and extend earlier studies on global bifurcations in the FitzHugh-Nagumo model, in particular those reported by Georgescu, Roc\c soreanu and  Giurgi\c teanu,
{\em Global Bifurcations in FitzHugh-Nagumo Model}, Trends in Mathematics: Bifurcations, Symmetry and Patterns (2003). \\ \\ 

\end{abstract}

\\


\section{Introduction}
Information in nerve fibers is encoded through \emph{action potentials} (electrical membrane changes).  A. Hodgkin and A. Huxley, in 1939, succeeded in recording an action potential from within a neuron for the first time, using microelectrodes inserted into the squid's giant axon. In their experiments, they demonstrated that the  action potentials arise from two effects:
\begin{itemize}
\item a fast inward current mediated by sodium (Na\(^+\)) ions, and
\item a slow developing outward current produced by potassium (K\(^+\)) ions.
\end{itemize}
 
Hodgkin and Huxley showed that the membrane’s permeabilities to $\text{Na}^+$ and $\text{K}^+$ are controlled separately, with their conductances varying in time and as functions of the membrane voltage. They formulated a quantitative description of the action potential in a single neuron, now referred to as the Hodgkin-Huxley (HH) model \cite{HoHu}.
 Details and an overview may be found in the review \cite{CPRG2024}.\\

\subsection{From the Hodgkin-Huxley model to the FitzHugh-Nagumo equations}
Hodgkin and Huxley represented the observed current variations using probabilistic variables that describe the closing and opening of ion channels. The Hodgkin-Huxley framework characterises the electrical dynamics of these membrane channels by modelling the flow of $\text{Na}^+$ and $\text{K}^+$ ions across the membrane.
In the HH framework, the dynamics are governed by four ODEs, each one associated with a specific state variable -- \cite[Eq. 1]{HoHu}.
This system is  nonlinear, making it analytically difficult to solve explicitly. 
Nonetheless, numerical simulations allow investigation of specific features and general dynamics, including the presence of \emph{excitability} and \emph{oscillations}.

The ability of the HH model to accurately reproduce these excitable dynamics represented a key insight, offering a 
 complete
explanation of the recorded action potential and the mechanisms underlying neural excitability (cf. \cite{HoHu}).
While the HH model successfully reproduces
many neuronal physiological phenomena, it is very complex.  In \cite{Fitzhugh}, Richard FitzHugh introduced a simplified model to capture the dynamics of neuronal excitability, which was  refined by Jinichi Nagumo  two years later \cite{Nagumo}. This system is known as the \emph{FitzHugh-Nagumo} (FHN) model.  
 FitzHugh concentrated the study on preserving the key dynamical characteristics of the HH model, specifically  and oscillatory and excitability behaviour.

\subsection{The FitzHugh-Nagumo equations}
FitzHugh started 
with the  oscillator equation (introduced in 1920) by Balthasar Van der Pol  \cite{vander}, 
which supports oscillations resembling relaxation-like dynamics.
This means there are periods of ``low'' and ``high'' states characterised 
through fast
transitions between them. 

The Van der Pol (VdP) equation is derived from the basic differential equation describing a damped harmonic oscillator. 
He replaced the constant damping by a  damping function 
that has a quadratic dependence on \(x,\) resulting in the nonlinear differential equation
(\(\dot{x}\) 
denotes
the first derivative of the physical position \(x\) with respect to \(t\) and \(c\in \RR^+ _0\)):
\begin{eqnarray*}
\ddot{x} +c(x^2-1) \dot{x} +x=0.
\end{eqnarray*}

The system exhibits ``effective" damping only when $|x|<1$; for \(x^2 > 1\), the nonlinear term leads to \emph{amplification}. To analyse the dynamics of the VdP equation, the Li\'enard map can be employed \cite{Linard}
\begin{eqnarray*}
y= \frac{\dot{x}}{c} +\left(\frac{x^3}{3}-x\right),
\end{eqnarray*}
giving rise to a system of two differential equations (on the plane):
\begin{eqnarray}
\label{FHN_intro}
\left\{ 
\begin{array}{l}
\dot{x}= \displaystyle c\left[y-\left(\frac{x^3}{3}-x\right)\right] \\ \\
 \dot{y}=  \displaystyle -\frac{1}{c}x,
\end{array}
\right.
\end{eqnarray}
from which the \emph{time-scale separation} of the equations becomes evident.

For $c\gg 1$, the first variable changes rapidly of the order of \(\mathcal{O}(c)\), whereas the second evolves much more slowly, \(\mathcal{O}(1/c),\) where \(\mathcal{O}\) denotes the standard \emph{Landau notation} (given a real map $f: \RR \to \RR$ and $c\in 
\RR$, we say that $f$ is of order $\mathcal{O}(c)$ if there exist $M\in \RR^+$ and $x_0\in \RR$ such that $|f(x)|<M c$, for  $x>x_0$). 
Based on
the VdP oscillator, FitzHugh considered the equation:

\begin{equation} \label{2d_intro}
  \left\{ 
\begin{array}{l}
\dot{x}=  \displaystyle c\left[y-\left(\frac{x^3}{3}-x\right)+z\right] \\ \\
 \dot{y}=  \displaystyle -\frac{1}{c}(x-a+by)
\end{array}
\right.
\end{equation}
where (see  \cite{CPRG2024} and Remark \ref{Rem1}):
\begin{enumerate}
\item   the parameter \(z\) mimics  the membrane current density;
\item the variable \(x\) 
corresponds to
the membrane voltage and the $\text{Na}^+$ activation and 
\item the variable \(y\) is related to the $\text{Na}^+$ inactivation and the $\text{K}^+$ activation. 
\end{enumerate}

Equations \eqref{2d_intro} act as an \emph{activator-inhibitor} framework.  
Arimoto, Nagumo, and Yoshizawa
proved that it can be realized as an electrical circuit \cite{Nagumo}.
Although the FHN model is primarily used to describe neuronal and cardiac systems, it can also be applied in a range of other biological contexts. 

The FHN equations provide a simplified framework for representing interconnected 
negative and positive
feedback loops, enabling the generation of a wide range of  responses such as switches, pulses, and stable oscillations \cite{Tsai}. 
Within the wide range of models used to
investigate cardiac cells, the FHN model is notable for being one of the simplest and most extensively researched, effectively representing the general dynamic characteristics of cardiac cells \cite{Alonso, Abid}. \\

\begin{rem}
\label{Rem1}
 System  \eqref{2d_intro} mimics  an excitable system and there are \textbf{limitations} regarding the interpretation of the original variables of the HH model.
In the rest of the paper, we will take a cautious approach and we stay within a mathematical motivation grounded in bifurcations and dynamical systems. \\
\end{rem}

 \subsection{Novelty}
  
Results in  \cite[pp. 180]{fitzhughbook} synthesised  the global bifurcation diagram for the FHN model. 
The result comes from integrating all the local bifurcation diagrams obtained in Chapters 2--4, as if fitting together a ``huge puzzle".
See  Section 2.4 of \cite{RS_proceedings}.
 
One of the novelties of this paper
is the dissection of 
the bifurcations diagrams of
\cite[pp. 180]{fitzhughbook}  and provide a more complete bifurcation analysis of model \eqref{2d_intro}, namely  the location of Hopf, Pitchfork and Double-zero bifurcations, in the cases that will be specified in Subsection~\ref{ss:structure}.   
 Besides  the case-study $c=\mathcal{O}(1)$, we also study the asymptotic case $$c \to \pm \infty \qquad ( \Leftrightarrow \quad \varepsilon = 1/c^2 \to 0),$$ where \emph{canards} are observed \cite{KS2001a, Kuehn}, and we finish the analysis with the  compactification of the phase portrait of \eqref{2d_intro} on the Poincar\'e disc. 

To the best of our knowledge, this procedure applied to equation \eqref{2d_intro} is  new,  providing additional information   about the trajectories which 
come from or tend to infinity.
We classify topologically 
all phase portraits of system \eqref{2d} ``near infinity'', represented on the Poincar\'e disc
(i.e. 
within the compactification of $\RR^2$ by including the circle $\mathbb{S}^1$ at infinity). 
Note that we  \textbf{do not} claim to have derived complete phase portraits. \\

\subsection{Structure} \label{ss:structure}
 This paper is structured as follows: Section \ref{sec2} introduces the definitions and preliminary concepts of 
quasi-homogeneous blow-up and Poincar\'e compactification. 
 In Section \ref{sec3}, we investigate the finite equilibria of \eqref{2d}  
  for the following three  cases: \\
\begin{description}
\item[Case A] \(a=0\)
\item[Case B] \(b=0\)
\item[Case C] \(a \neq 0,\, 0<b<1\) \\
 \end{description}
and investigate possible   bifurcations, including Pitchfork, Double Homoclinic and  Hopf bifurcations. The general bifurcation analysis of \eqref{2d_intro}  is difficult to tackle; however we have been able 
 to make some  progress in \textbf{Cases A, B} and \textbf{C.}
In Section \ref{s:canards}, we analyse the asymptotic dynamics  of  \eqref{2d_intro}   when $c \to \pm \infty$ and we relate the existence of \emph{canards} with the periodic solution emerging from the Hopf bifurcation.

Section \ref{sec4} is dedicated to studying the dynamics ``near infinity'' for \eqref{2d_intro} on the Poincar\'e disc.
The phase portraits for each connected component in  \textbf{Cases A, B} and \textbf{C} are 
 illustrated.  This complements the work started in \cite{GlobalFN}. Section \ref{s: conclusion} finishes the paper. \\

In this paper, we aim to provide a self-contained exposition, integrating all topics pertinent to the proofs. Illustrative figures have been included for clarity, and all results are supported by numerical simulations performed using \emph{Matlab R2015b} and \emph{pplane9}.\\

\section{Preparatory section} \label{sec2}
 
This section provides definitions for polynomial vector fields on \(\RR^2\) that will be used
in the sections that follow.
Let \(f_{(a,b,c)}\) be a smooth vector field on \(\RR^2\) 
whose flow is given by
the unique solution \(X(t)=\varphi(t, X)\in \RR^2\) of the three-parameter family
 
\begin{equation}\label{general}
\dot{X}=f_{(a, b, c)} (X), \qquad X(0)=X_0 \in \RR^2,
\end{equation}
where
 $X=(x,y)\in \RR^2$   and 
\((a, b, c)\in \RR \times \RR \times \RR\backslash \{0\}\).  \\

\subsection{Useful terminology}
 
The \emph{center manifold} of a non-hyperbolic equilibrium $p\in \RR^2$ of \eqref{general} is the set of solutions whose behaviour around $p$ is controlled neither by the exponential attraction of the stable manifold nor by the exponential repulsion of the unstable manifold. 

 If  \(Df_{(a, b, c)}\) evaluated at an equilibrium $p$,  \(Df_{(a, b, c)}(p)\),  has an eigenvalue with zero  real part, the center manifold plays an important role and this is the set where \emph{bifurcations} might occur.   

For a  triple \((a, b, c)\in \RR \times \RR \times \RR\backslash \{0\}\), an equilibrium point \(p\) of \eqref{general} is called \emph{nilpotent singularity} when all eigenvalues of Jacobian matrix \(Df_{(a,b,c)}(p)\) are zero but \(Df_{(a,b,c)}(p)\neq 0\), and is called \emph{linearly zero} when \(D f_{(a,b,c)}(p)\equiv 0.\) 
Further, a singular point \(p\) is termed \emph{semi-hyperbolic} if precisely one of eigenvalues of \(Df_{(a,b,c)}(p)\) is equal to \(0.\)

Throughout this paper, we study the Double-zero (DZ) \emph{singularity}   with symmetry \(\mathbb{Z}_2\) for a family of differential equations  corresponding to the case  described in  \cite[pp. 400]{wiggins}.  The unfolding of this \emph{singularity} of codimension two involves lines of  \emph{Belyakov transitions}, \emph{Pitchfork}, \emph{Hopf}, \emph{Double Homoclinic
and saddle-node   bifurcation of limit cycles}. 

 A Hopf bifurcation is termed \emph{supercritical} when the first Lyapunov coefficient of its normal form is negative, leading to the emergence of a stable small-amplitude limit cycle, and \emph{subcritical} when this coefficient is positive, in which  an unstable small-amplitude limit cycle is created \cite{Verduzco}.
 For full details on these bifurcations and the criteria that guarantee their appearance, see \cite{Kuznetsov, CR2023, wiggins}.\\

\begin{rem}
In this paper, we refer to an equilibrium point $p$ as undergoing a \textit{Belyakov transition} if at least one pair of eigenvalues of $Df_{(a,b,c)}(p)$  changes from real to complex conjugate or vice versa, while the sign of their real part remains unchanged. Although such a transition is typically considered in higher-dimensional systems, we adopt this terminology here for the sake of brevity and clarity. \\
\end{rem}

\subsection{Quasi-homogeneous blow-up method}\label{sec:vertical}
\emph{Quasi-homogeneous blow-ups} are used to analyze and simplify the behaviour of differential equations near non-hyperbolic singularities by rescaling the variables according to their natural weights. This technique transforms the singularity into a more regular, analyzable object. Although homogeneous blow-ups
 \begin{equation*}
(\theta, \rho) \mapsto \left( \rho \cos \theta,\, \rho \sin \theta \right),
\end{equation*}
suffice for the study of isolated singularities of an analytic vector field, the \emph{quasi-homogeneous} approach  (scaling the variables with different degrees) is considerably more efficient. Here we review this method briefly -- for more details and examples, see \cite[Section 3.3]{llibrebook}.\\

Consider the polynomial vector field $f=(P, Q)$ associated with the planar differential equation: 
\begin{equation}\label{rs}
\dot{x}_1=P(x_1,x_2), \qquad
\dot{x}_2=Q(x_1,x_2).
\end{equation}
In the language of \emph{singularity theory}, we have:
\begin{equation*}
f(x_1,x_2)= P(x_1, x_2) \frac{\partial}{\partial x_1} + Q(x_1, x_2) \frac{\partial}{\partial x_2},
\end{equation*}
  with an isolated singular point at the origin $(0,0)$. Since $P$ and $Q$ are polynomial we may find $p_{ij}, q_{ij} \in \RR$ such that:
\begin{equation}
\label{P and Q}
P(x_1, x_2) = \sum_{i + j \geq 1} p_{ij} x_1^i x_2^j \qquad\text{and}\qquad Q(x_1, x_2) = \sum_{i + j \geq 1} q_{ij} x_1^i x_2^j.
\end{equation}

\bigbreak

\begin{defn}
The \emph{quasi-homogeneous blow-up} is defined by 
\begin{equation*}
\phi : \mathbb{S}^1 \times \mathbb{R} \to \mathbb{R}^2, \quad 
(\theta, \rho) \mapsto \left( \rho^{\alpha} \cos \theta,\, \rho^{\beta} \sin \theta \right)=(x_1,x_2),
\end{equation*}
  where
\( \alpha\) and  \(\beta\)
are natural numbers. 
\end{defn}
Setting
\(\alpha=\beta=1,\) recovers the \emph{standard homogeneous} case. 
Doing  the quasi-homogeneous blow-up (of the previous definition) on 
$f=(P, Q)$, it yields the following differential equation:  
\[
\dot{\rho} = \frac{P \rho^{\beta} \cos\theta  + Q \rho^{\alpha} \sin\theta }{\rho^{\alpha+\beta-1} (\alpha \cos^2\theta + \beta \sin^2\theta)}, \qquad
\dot{\theta} = \frac{\alpha Q \rho^{\alpha}  \cos\theta - \beta P \rho^{\beta}  \sin\theta}{\rho^{\alpha+\beta} (\alpha \cos^2\theta + \beta \sin^2\theta)},
\]
in which
\(P\) and \(Q\) are 
computed
at \((\rho^{\alpha} \cos\theta, \rho^{\beta} \sin\theta)\). A common factor \(\rho^d>0\) can be cancelled for a suitable choice of \(d>0\). \\

Throughout this article, we will be interested in  \emph{quasi-homogeneous directional blow-ups in a given direction}. This technique isolates and desingularizes the dynamics of trajectories approaching the singularity along that direction, creating a chart where the vector field becomes simpler and its behaviour can be analyzed separately from other directions. \\
 
\begin{defn}
The \emph{quasi-homogeneous directional blow-up 
along
the positive} (resp. negative) \(x\) 
axis
is given by the map:
\begin{align*} 
&\overline{f }^{x}_{+}: \quad (\bar{x}, \bar{y}) \mapsto (\bar{x}^{\alpha}, \bar{x}^{\beta} \bar{y}), \\
&\overline{f }^{x}_{-}: \quad (\bar{x}, \bar{y}) \mapsto (-\bar{x}^{\alpha}, \bar{x}^{\beta} \bar{y}), 
\end{align*}
and the \emph{quasi-homogeneous directional blow-up 
along
the positive} (resp. negative) \(y\) 
axis
is described by:
\begin{align*} 
&\overline{f }^{y}_{+}: \quad (\bar{x}, \bar{y}) \mapsto (\bar{x} \bar{y}^{\alpha}, \bar{y}^{\beta}),\\
& \overline{f }^{y}_{-}: \quad (\bar{x}, \bar{y}) \mapsto (\bar{x} \bar{y}^{\alpha}, -\bar{y}^{\beta}),
\end{align*}
where \(\bar{x}\) and \(\bar{y}\) denote new variables and   \( \alpha, \beta \in \mathbb{N}.\) 
\end{defn}
\bigbreak
 
Among all possible pairs   $(\alpha, \beta)\in \NN\times \NN$ of the definition of quasi-homogeneous blow-ups, a natural question is: 
 which choice is the most ``efficient'' (in the sense that we desingularize the equilibria in the minimum number of steps) for performing a quasi-homogeneous blow-up to simplify and resolve a singularity? 
 
 Selecting such an efficient pair balances the relative importance of the variables, making the leading-order dynamics near the singularity clear. The \emph{Newton polygon} provides a systematic tool to guide this choice, which we now introduce (adapted from \cite{llibrebook}). \\
\bigbreak
\begin{defn}
Under the notation of \eqref{P and Q}, we define the \textit{support} of $f$ as the subset
\begin{equation*}
S = \{ (m - 1, n) \mid p_{mn} \neq 0 \} \cup \{ (m, n - 1) \mid q_{mn} \neq 0 \} \subset \mathbb{R}^2,
\end{equation*}
where the point \((-1, n)\) corresponds to the monomial \(p_{0n} x_2^n\) and the point \((m,-1)\) corresponds to the monomial \(q_{m0}x_1^m.\) Furthermore,  \((0, 0)\) is associated with the monomials \(p_{10}x_1\)
and \(q_{01}x_2.\)
\end{defn}
 \bigbreak 
 \begin{defn}
The \textit{Newton polygon} associated to $f$ is given by the convex hull 
 \begin{equation*}
 \bigcup_{(r,s) \in S} \{ ({r}^*, {s}^*) \in \mathbb{R}^2 \mid r\leq {r}^*,\, s\leq {s}^*  \}.
\end{equation*}
\end{defn}
The \textit{Newton polygon} associated to $f$  consists of the union \( \gamma \) of the points/segments  \( \gamma_k \), 
enumerated in order from left to right.
For each 
$1\leq k\leq n$,
the segment $\gamma_k$ lies on the straight line
\(
\alpha_k r + \beta_k s = d_k
\)
where $\alpha_k$ and $\beta_k$ are coprime integers and $d_k$ is a constant.
Since $f$ has an isolated singularity at the origin then, at least one of the two points $(0,s)$ or $(-1,s)$,
belongs to the support $S$
and similarly, 
 at least one of two points $(r,-1)$ or $(r,0)$
lies in $S$, for some $r$. 
Thus a first segment $\gamma_1$ always exists in the Newton diagram ($\Leftrightarrow$ it is non-empty). \\

The Newton diagram identifies the dominant monomials near a singularity; the slopes of its compact edges determine the weight ratio $\beta/\alpha$, $\alpha \neq 0$, used to choose the exponents of a quasi-homogeneous blow-up. \\

\begin{prop}\cite[Section 3.3, adapted]{llibrebook}.
Under the previous conditions, the \emph{most efficient} pair $(\alpha, \beta)$ for the quasi-homogeneous blow-up  corresponds to the edge of the Newton polygon with minimal slope (in absolute value). \\
\end{prop}

\subsection{Poincar\'e compactification} \label{ss:compactification(prel)}

Here, following Chapter 5 of \cite{llibrebook}, we present a concise overview of the Poincar\'e compactification 
to study how trajectories of a planar differential system behave ``near infinity''. This technique  will be used in Section \ref{sec4} of the present paper. 

To construct the phase portrait,
one would ideally need to consider 
the complete real plane \(\mathbb{R}^2,\) 
 that is often impractical. However, if the vector field is defined by polynomial functions, Poincar\'e compactification can be applied. This method allows the phase portrait to be represented within a finite region and also provides control over orbits that approach or come from infinity. 
 Consider the polynomial vector field \(f=(P,Q)\),
where \(d_1\) and \(d_2\) are the algebraic degrees of \(P\) and \(Q\) respectively and  \(d={\max} \{d_1,d_2\}.\)  \\

\subsubsection{Preliminaries for the construction}
 
First identify \(\RR^2\) with the plane \(\Pi\) in \(\RR^3\) defined 
as  \((y_1, y_2, y_3)=(x_1, x_2, 1).\) 
The sphere 
$$ \mathbb{S}^2=\left\{(y_1, y_2, y_3)\in \RR^3: \sum_{i=1}^3 {y_i}^2=1\right\}$$
is called the \emph{Poincar\'e sphere} and  is tangent to 
 \(\Pi\) at \((0,0,1).\)
The sphere may be written as $$H^+\cup H^-\cup \mathbb{S}^1$$ where
\begin{equation*}
H^+=\{(y_1, y_2, y_3)\in \mathbb{S}^2:y_3>0\} , \quad H^-=\{(y_1, y_2, y_3)\in \mathbb{S}^2:y_3<0\} \end{equation*}
 and the equator $\mathbb{S}^1$ given by:
\begin{equation*} 
\mathbb{S}^1 = \{(y_1, y_2, y_3)\in \mathbb{S}^2:y_3=0\}.
\end{equation*}
If \(x=(x_1, x_2)\in \RR^2\) define \(\Delta(x)\) as   \(\sqrt{x_1^2+x_2^2+1}\neq 0 .\) \\

\begin{defn}
The projection of \(\mathbf{X}: \RR^2\to \mathbb{S}^2 \)
is defined with
the central projections
\begin{equation*} 
f^{\pm}: \RR^2\to \mathbb{S}^2 
\end{equation*}
where \(f^{\pm}(x_1, x_2)\) is the intersection point of the line that passes through
\((0,0)\)
and \((x_1, x_2)\in \RR^2\), with $H^{\pm}$. More precisely:
\begin{equation*} 
f^+(x_1, x_2)= \left(\frac{x_1}{ \Delta(x)},\frac{x_2}{ \Delta(x)}, \frac{1}{ \Delta(x)}\right).
\end{equation*}
Analogously, we define \(f^-\) (substituting $H^+$ by $H^-$) as:
\begin{equation*} 
f^-(x_1, x_2)= \left(-\frac{x_1}{ \Delta(x)},-\frac{x_2}{ \Delta(x)},-\frac{1}{ \Delta(x)}\right).
\end{equation*}
\end{defn}

For $x=(x_1, x_2)\in \RR^2$, the induced vector
field on 
 $H^{\pm}$
is defined, respectively, by: 
 
 \begin{align*}
\overline{ \mathbf{X}}(f^{\pm}(x)) &= Df^{\pm}(x)  \mathbf{X}(x).\\
 \end{align*}
 

The infinite points of \(\RR^2\)
(each direction is associated with two points) are in
bijective 
 associated to
the points of the equator of \(\mathbb{S}^2.\) 
In the next subsection, we are going 
 to generalize
the induced vector field 
 \(\overline{ \mathbf{X}}: \mathbb{S}^2\backslash \mathbb{S}^1\rightarrow \mathbb{S}^2.\)\\
\begin{defn}
The extension of \(\mathbf{X}\) to \(\mathbb{S}^2\) is referred to as its \emph{Poincar\'e compactification}, denoted by \(p(\mathbf{X}).\)  \\
 \end{defn}
  
\subsubsection{About the construction} We use smooth charts to make calculations. For \(y=(y_1, y_2, y_3) \in \mathbb{S}^2,\) we 
utilize
 local charts:
\begin{eqnarray*}
U_l =\{y\in \mathbb{S}^2 : 0<y_l \}, \quad V_l =\{y\in \mathbb{S}^2 : 0>y_l\}, \quad l=1,2,3. 
 \end{eqnarray*}
The 
associated
 local maps $\Phi_l : U_l \to \RR^2$ and $\Psi_l : V_l \to \RR^2$ 
 are described as
\begin{eqnarray*}
 \Phi_l(y) = -\Psi_l(y) = \left(\frac{y_m}{y_l},\frac{y_n}{y_l}\right), \quad  \quad m < n \quad \text{where} \quad m,n \neq l.
 \end{eqnarray*}
Consider
$z = (u,v)$ as the value of $\Phi_l(y)$ or $\Psi_l(y)$ for any $l\in \{1,2,3\}$,
such that $(u,v)$  plays different roles depending on the local chart we are considering.
In any chart, points on \(\mathbb{S}^1\) correspond to those with \(v=0\).\\

\subsubsection{Explicit expressions}

In the local chart $(U_1,\Phi_1),$ $p(\mathbf{X})$ is expressed as
\begin{equation}\label{comp1}
\dot{u}=\left[
Q\left(\frac{1}{v},\frac{u}{v}\right)-uP\left(\frac{1}{v},\frac{u}{v}\right)\right] v^d, \hspace{0.5cm}
\dot{v}=- P\left(\frac{1}{v},\frac{u}{v}\right)v^{d+1},
\end{equation} 
 for \((U_2, \Phi_2)\) 
 :
\begin{equation}\label{comp2}
\dot{u}=\left[P\left(\frac{u}{v},\frac{1}{v}\right)-
u Q\left(\frac{u}{v}, \frac{1}{v}\right)\right]v^d, \hspace{0.5cm}
\dot{v}= -Q\left(\frac{u}{v},\frac{1}{v}\right)v^{d+1},
\end{equation}
and 
for \((U_3, \Phi_3)\) is:
\begin{equation}\label{comp3}
\dot{u}=P(u,v), \hspace{0.5cm} \dot{v}=Q(u,v).
\end{equation}
 
 \begin{rem}
For \(i=1,2,3\), the 
representations
for  \(p(\mathbf{X})\) in 
the local chart \((V_i,\Psi_i)\) are identical to those in 
\((U_i,\Phi_i),\) up to a multiplication  by \(({ -1})^{d-1}\).
\end{rem}

 To investigate \(\mathbf{X}\) throughout the complete plane \(\mathbb{R}^2\), 
 as well as the behaviour ``near infinity'', 
it is enough to
examine it on \(H^{+} \cup \mathbb{S}^1\), a region known as the \emph{Poincar\'e disc}.
All calculations can be  performed in the three charts 
$(U_i,\Phi_i)$ for \(i=1,2,3\)
whose expressions are given by
the formulas \eqref{comp1}, \eqref{comp2} and \eqref{comp3}.     \\
 \begin{defn}
 We refer to singular points of $\mathbf{X}$ as \emph{finite} when the corresponding points of $p(\mathbf{X})$ are in \(\mathbb{S}^2 \setminus \mathbb{S}^1,\) and as \emph{infinite} when they belong to \(\mathbb{S}^1\). \\
\end{defn}

\subsubsection{Useful consequences}
Following Chapter 5 of 
\cite{llibrebook}, we list some useful consequences of the theory described above. \\
\begin{enumerate}
\item 
Whenever $y \in \mathbb{S}^1$ is an infinite singular point, 
so is $-y$.
\item  
The local behaviour near $-y$ is that near $y$ scaled by $(-1)^{d-1}$.
\item 
 Because infinite singular points 
 arise in
pairs, it is sufficient to examine just half of them. The rest can be inferred via the vector field's degree.
\item 
 The integral curves on \(\mathbb{S}^2\) are symmetric around \((0,0)\), so it suffices to depict the flow of \(p(\mathbf{X})\) in the closed northern hemisphere.\\
\end{enumerate}
 
\subsubsection{Infinite singular points}
The theory described here will be useful in Section  \ref{sec4} of the present paper.
Our goal is to analyse the local phase portrait near infinite singular points. To this end, we select a singular point \((u,0)\) and begin with examining the linear part of the vector field \(p(\mathbf{X})\).
 Let $P_i$ and $Q_i$ be homogeneous monomials of degree $i \in \mathbb{N}_0$,
 $i = 0,1,...,d$, such that \eqref{P and Q} may be written as:
  \begin{equation*} 
 \left\{ 
\begin{array}{l}
P= \displaystyle \sum_{i=0}^{d_1} P_i,\\\\
Q=\displaystyle \sum_{i=0}^{d_2} Q_i.
 \end{array}
\right.
\end{equation*}
 Hence, if $(u,0) \in  (U_1 \cup V_1) \cap \mathbb{S}^1$ is an infinite singular point 
  for 
$p(\mathbf{X})$ then:
 \begin{eqnarray*}
 F (u) \equiv   - uP_d (1, u)+Q_d (1, u) = 0.
\end{eqnarray*}
Analogously, if 
$(u,0) \in (U_2 \cup V_2) \cap \mathbb{S}^1$ is an infinite singular point for $p(\mathbf{X})$ 
then:
\begin{eqnarray*}
 G (u) \equiv   - uQ_d ( u,1)+P_d ( u,1) = 0.
\end{eqnarray*}
The reverse of the previous assertions are also true. The Jacobian of the vector field $p( \mathbf{X})$ at 
$(u,0)$ is
\begin{eqnarray*}
 \left(\begin{array}{cc}
F'(u) & Q_{d-1}(1,u)-uP_{d-1}(1,u) \\ 
\\
0 &  -P_d(1,u)
\end{array}\right)  
\quad \text{or}\quad 
 \left(\begin{array}{cc}
G'(u) & P_{d-1}(u,1)-uQ_{d-1}(u,1) \\ 
\\
0 &  -Q_d(u,1)
\end{array}\right)  
\end{eqnarray*}
if $(u, 0)$ 
lies in
$U_1 \cup V_1$ or $U_2 \cup V_2$, respectively. \\

Our discussion is concentrated on \emph{isolated singularities} in the equator.
According to \cite{llibrebook}, 
only \emph{nodes} and \emph{saddles} can arise as hyperbolic singularities at infinity.
All semi-hyperbolic singular points may also 
emerge
at infinity. 
In the case that 
one of these 
singular points
at infinity is 
  a saddle (topologically),
then 
 the line
\(v=0,\) corresponding to the equator of \(\mathbb{S}^2\), serves as a stable, unstable, or center manifold. This property also applies to semi-hyperbolic points of  saddle-node type, whose hyperbolic sectors may  be 
divided into two
 distinct forms
 depending on 
the Jacobian matrix written in the charts $U_1$ or \(U_2.\)
These matrices may be either
\begin{eqnarray*} 
 \left(\begin{array}{cc}
0 & \star\\ 
\\
0 &a
\end{array}\right)  
\quad \text{or}\quad 
 \left(\begin{array}{cc}
a & \star \\ 
\\
0 & 0
\end{array}\right) 
\end{eqnarray*}
where
$a \neq 0$ and $\star \in \RR$.  
Nilpotent equilibria exhibit behaviour at infinity that is considerably more complex than those of hyperbolic equilibria. 
Blow-up techniques of Subsection \ref{sec:vertical} are needed to study them.

\section{Finite equilibria and bifurcation analysis}\label{sec3}
Our object of study is the analysis of the  family of differential equations \eqref{2d_intro}, which can be written as:
\begin{equation} 
 \label{2d}
  \left\{ 
\begin{array}{l}
\dot{x}= \displaystyle c\left[y-\left(\frac{x^3}{3}-x\right)\right] \\ \\
\dot{y}= \displaystyle -\frac{1}{c}(x+by-a)
\end{array}
\right.
\end{equation}
where \(a, b\in \RR\) and \(c\in \RR\backslash\{0\}.\)
 We can absorb the parameter \(z\) of \eqref{2d_intro} into the variable through the change of variable
\(y \to y+z\) and \(a\to a-bz.\) The model \eqref{2d} is invariant under the transformation \((t,c)\to (-t, -c)\) so we can restrict to the case \(c>0\), which means that the stability of the equilibria for \(c>0\) is the reverse of that of \(c<0\). \\

Let us denote by \(f_{(a,b,c)}: \mathbb{R}^2\to \mathbb{R}^2\) the vector field associated with \eqref{2d}. 
For \(a, b\in \RR\) and $c\in \RR\backslash\{0\}$, the Jacobian matrix of \(f_{(a,b, c)}\) at a general point \((x,y)\in \RR^2\) is given by:
\begin{eqnarray}
 Df_{(a,b, c)} (x,y)=
\left(\begin{array}{cc}
c-cx^2 & c \\ 
\\
-1/c & -b/c
\end{array}\right).  \label{JE2_DZb}
\end{eqnarray}
 It is immediate to deduce:
\begin{lem}
\label{lemma3.1}
The divergence of  (\ref{2d}) is  given by  \(c-cx^2-b/c,\) which is strictly negative for \(b>c^2.\)
\end{lem}

 By the Bendixson criterion   \cite[Theorem 7.10]{llibrebook}, the negativity of the previous result ensures that the system (\ref{2d}) cannot have periodic orbits and, consequently, no limit cycles in the open set defined by the inequality \(b>c^2.\) From now on,  we divide the analysis into Cases \textbf{A}, \textbf{B} and \textbf{C} (see Subsection \ref{ss:structure}), depending on the parameters $a,b$ and $c$ of (\ref{2d}). \\

\subsection{Case A (\(a=0\))}
Consider the invertible linear map in \(\mathbb{R}^2\) defined by \(\kappa(x,y)=(-x, -y)\), whose action on \(\RR^2\) is isomorphic to that of $\mathbb{Z}_2$ (rotation of $\pi$ around the origin). Since $$f_{(0,b,c)}\circ \kappa = \kappa \circ f_{(0,b,c)},$$ we may say that: 
 
\begin{lem}
The vector field \(f_{(0, b,c)}\) is \(\mathbb{Z}_2(\kappa)\)--equivariant.
\end{lem}

For all $a=0$, $b  \in \RR$ and \(c\in \RR\backslash\{0\},\) we know that \( E_1= (0,0) \) is an equilibrium of  \eqref{2d}. If \(b\in (-\infty, 0)\cup [1, +\infty),\) then system  \eqref{2d} has two extra equilibria given explicitly by  
\begin{eqnarray}\label{E2E3_def}
&E_2=\kappa(E_3)= \left(-\sqrt{\dfrac{3(b-1)}{b}},\dfrac{1}{b}\sqrt{\dfrac{3(b-1)}{b}} \right)\\ \quad \text{and}\quad
\nonumber &E_3=\kappa(E_2)= \left(\sqrt{\dfrac{3(b-1)}{b}},-\dfrac{1}{b}\sqrt{\dfrac{3(b-1)}{b}} \right).
\end{eqnarray}
From  \eqref{JE2_DZb},  we deduce that the trace and determinant operators of \(Df_{(0,b, c)},\)  evaluated at \(E_i,\) \(i=1,2,3,\) are given and denoted, respectively, by:
\begin{align}
 \nonumber \mbox{Tr}(E_1)&=c-\dfrac{b}{c}, \qquad\qquad\qquad\qquad \mbox{Det}(E_1)=1-b,\\\nonumber
\mbox{Tr}(E_{2,3})&=-\dfrac{2 c^2 b-3 c^2+b^2}{c b}, \qquad \mbox{Det}(E_{2,3})=2(b-1).
\end{align}
It can be easily shown that $E_1$ is a hyperbolic saddle
for \(b>1\) (since \(\mbox{Det}(E_1)<0\)).
Using \eqref{JE2_DZb}, the jacobian matrix of the vector field \(f_{(0,1,\pm 1)}\)   at \(E_1\equiv E_2\equiv E_3\) is given by: 
\begin{eqnarray}\label{JE2_DZb2}
 &&\left(\begin{array}{cc}
\pm 1  & \pm 1 \\ 
\\
\mp 1 & \mp1
\end{array}\right) .
\end{eqnarray}
It is straightforward to see that
the matrix \eqref{JE2_DZb2} is non-hyperbolic and has a double zero eigenvalue.  Our main result relies on the existence of a \emph{Double-zero bifurcation} (DZ) with \(\mathbb{Z}_2(\kappa)\)--symmetry of $f_{(0, 1, \pm1)}$ at the singularity $E_1$.
 
\begin{rem}
Consider the following linear change of coordinates, \(\Psi,\) and its inverse:
\begin{eqnarray*}
\Psi (x,y)= (x+y,-x) \quad \text{and} \quad \Psi^{-1}(x,y)=(-y, x+y).
\end{eqnarray*}
One may easily check that: 
\begin{eqnarray*}
 \Psi^{-1}\circ f_{(0,1,\pm 1)}\circ \Psi(x,y)&=& \Psi^{-1}\circ f_{(0,1,\pm 1)}(x+y, -x)\\
 &=& \Psi^{-1}\left(-\frac{(x\pm y)^3}{3}+y, -y\right)\\
 &=& \left(y, -\frac{1}{3} (x\pm y) ^3\right),
\end{eqnarray*}
 \end{rem}
so that the equations in Jordan form 
are represented by
\begin{equation}
\label{NFDZ} 
   \left\{ 
\begin{array}{l}
\dot{x}=y  \\ \\
\dot{y}=-\frac{1}{3} (x\pm y) ^3,
\end{array}
\right.
\end{equation}
  which contains part of the truncated form of degree 3 of the versal deformation of a  \emph{ Double-zero bifurcation} with \(\mathbb{Z}_2 \)--symmetry  \cite[pp. 400, Case $c=-1$]{wiggins}.  The next result deals with the description of all bifurcation curves passing through the bifurcation point $(b,c)=(1, \pm 1)$ and that characterises a  DZ bifurcation  of \(f_{(0,b,c)}\) with \(\mathbb{Z}_2 (\kappa)\)--symmetry. 

 We focus the proof on the
 analysis around the bifurcation parameter  \((b, c) =\left(1,  1\right)\); the analysis close to  \((b, c) =\left(1,  -1\right)\)   has a similar treatment (reversing time).\\
 
\begin{figure}[t!]
\begin{center}
{\includegraphics[width=.48\columnwidth,height=.45\columnwidth]{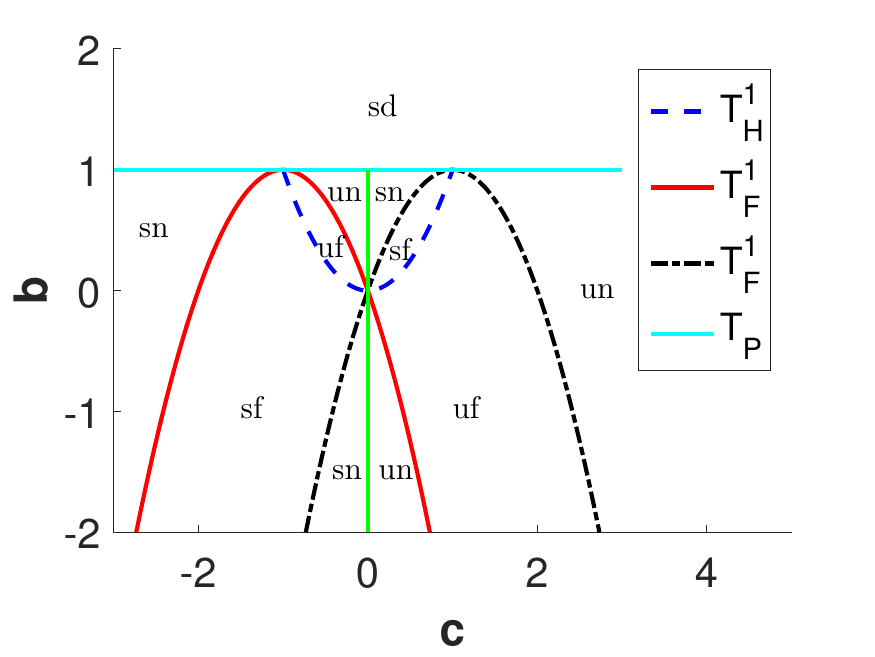}}
  {\includegraphics[width=.48\columnwidth,height=.45\columnwidth]{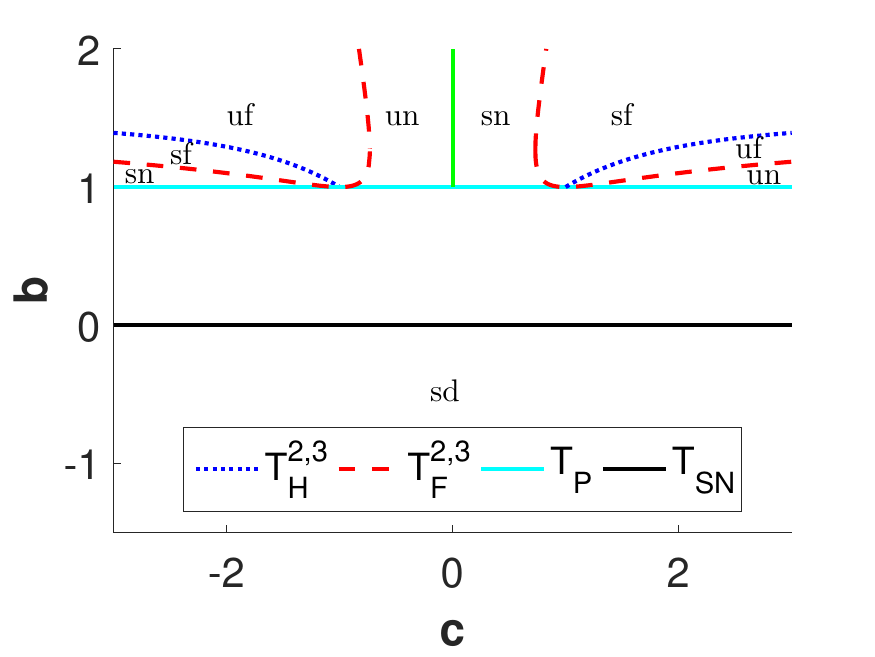}} 
 \end{center}
\caption{ \label{1fig} 
 \small   Illustration of Theorem \ref{th:main1}.  \textbf{Left:} The bifurcation curves associated to  \eqref{2d} for \textbf{Case A}.  Linear stability of \(E_1\) for \(f_{(0,b,c)}\). \textbf{Right:} Linear stability of \(E_2\) and \(E_3\) for \(f_{(0,b,c)}\).   The notation $T_P$, $T_H^1$, $T_F^1$, $T_P$  $T_H^{2,3}$, $T_F^{2,3}$  follows from Theorem \ref{th:main1} and \textbf{sd}, \textbf{uf}, \textbf{un}, \textbf{sf}, and \textbf{sn} refer to saddle, unstable focus, unstable node, stable focus and stable node.  
 }
\end{figure}

\begin{thm}
\label{th:main1}
Regarding to
the vector field \(f_{(0,b,c)},\) the equilibrium \(E_1\equiv E_2 \equiv E_3\) undergoes a  DZ bifurcation with \(\mathbb{Z}_2(\kappa)\)--symmetry at \((b, c)=(1,   1).\) The global representations of the transition/bifurcation curves in the 
parameters space
\((b,c)\in \RR \times \RR^+\) are as follows (schematic lines for \(a=0\) have been plotted  in Figures \ref{1fig} and \ref{CASEA}):\\ 
\begin{enumerate}
 \item The equilibrium \(E_1\) undergoes  a supercritical (with respect to $b$)  Pitchfork bifurcation   along the line:
 \begin{align}
 \nonumber
T_{P}=\{(b,c): b=1\}
\end{align} 
giving rise to the equilibria \(E_2\) and \(E_3\) defined in (\ref{E2E3_def}) for $b >1$.    \\ \\

\item The equilibrium \(E_1\) undergoes a Belyakov transition along the union of parabola: 
 \begin{align}
\nonumber
T^1_F=\{(b,c): b=-c^2\pm 2c\}.
\end{align}
 \\
 \item The equilibrium \(E_1\) undergoes a 
 supercritical
(with respect to $b$) Hopf bifurcation along part of the parabola:
 \begin{align*}
T_{H}^1=\left\{(b,c): b=c^2,\, |c|<1\right\}.
\end{align*}
\\

\item The equilibria \(E_2\) and \(E_3\) undergo a Belyakov transition along the algebraic curves: 
 \begin{align*}
T^{2,3}_F=\left\{(b,c): \left(b^2+2cb-2c^2b+3c^2\right)\left(b^2-2cb-2c^2b+3c^2\right)=0\right\}.
\end{align*}
\\
\item The equilibria
\(E_2\) and \(E_3\) undergo a 
subcritical
(with respect to $b$) Hopf bifurcation along the algebraic curves:
\begin{align*}
T^{2,3}_H=\left\{(b,c): b=c\left(-c\pm\sqrt{c^2+3}\right) \mbox{ and } b>1\right\}.
\end{align*}
  \\
 \item There is a  double homoclinic cycle to \(E_1\)  along the approximated curve:
\begin{eqnarray*}
DH:\quad 0 =   \frac{7b^2 + 10bc^2 - 17c^2}{15c^3} \sqrt{\frac{-7b^2 + 5bc^2 + 2c^2}{5b}}, \quad b \in (0, c). 
\end{eqnarray*}
 \\
\end{enumerate}
\end{thm}

\begin{rem}
  For every  pair \((b,c)\in \RR \times \RR^+\), the origin is always an equilibrium point.
When $b<0$, two additional nontrivial equilibria appear; they move away from the origin and diverge to infinity as  $b \to 0^{-}$.  For $0<b<1$, the origin is the only equilibrium. 
 Thus, although the number of equilibria changes as the line $b=0$ is crossed, this change results from equilibria escaping to (or emerging from) infinity rather than from a local bifurcation at the origin.\\
\end{rem}

\begin{rem}
\label{rem:SNL}
According to \cite[Theorem 3.5.1]{fitzhughbook} one knows that there is a \emph{Bautin bifurcation} along the curve parametrized by $c$ defined by 
\begin{equation}
\label{SNL}
 \left\{(a,b, c)\in \RR^2\times \RR^+: \left( \pm \frac{4}{3} \left( c\sqrt{c^2-1}-c^2+1\right) \sqrt[4]{1-\frac{1}{c^2}}, c^2 \pm  c\sqrt{c^2-1}, c\right), c\geq1 \right\}
\end{equation}

In particular, in the three-dimensional bifurcation space $(a, b, c)$, this gives rise to   a saddle-node bifurcation surface   associated with two non-hyperbolic cycles, generating the line SNL defined by $b=c^2\pm c\sqrt{c^2-1}$ and $a=0$  of Figure \ref{summary1}. The formula of  \cite{fitzhughbook} is explicitly written for $c\geq1$, however due to the invariance under the reversible transformation \((t,c)\to (-t, -c)\), it also holds for $c\leq -1$. 


The curve SNL of Figure \ref{summary1} does not make part of a generic unfolding of a DZ bifurcation with symmetry; this is why it is represented as a dashed line.

When $a=0$, the points $\textbf{Q}_{17}$ and $\textbf{Q}_{18}$ of \cite[pp. 180]{fitzhughbook} coincide with the point $(b,c)=(1,1)$. This coincidence yields a codimension-three DZ--bifurcation  with \(\mathbb{Z}_2(\kappa)\)--symmetry at $(b, c)=(1,   1)$, whose analysis is beyond the scope of the present paper.   A complete understanding of this bifurcation remains an open problem and is briefly discussed in Section \ref{s: conclusion}. \\
\end{rem}

\begin{rem}
Theorem \ref{th:main1} generalizes the analysis of \cite{GlobalFN} and 
\cite{fitzhughbook} by locating a codimension 2 bifurcation  with \(\mathbb{Z}_2(\kappa)\)--symmetry  at \((b,c)=(1,1)\) for the vector field \(f_{(0,b,c)}.\) For the sake of completeness,  
we present
 local numerical phase portraits on Table \ref{QA} of Appendix. All of them agree well with the theory.
 \\
\end{rem}
\bigbreak
\begin{proof}[Proof of Theorem  \ref{th:main1}]
The first statement of the result comes from the normal form \eqref{NFDZ}.
The eigenvalues of the Jacobian matrix \eqref{JE2_DZb} at   \(E_1,\) \(E_2\) and \(E_3,\) when the equilibria exist, are: \\
\begin{align*}
E_1:\quad &\lambda_{\pm}^1=\dfrac{c^2-b\pm\sqrt{(c^2+b)^2-4 c^2}}{2c},\\\nonumber
E_2, E_3:\quad &\lambda_{\pm}^{2,3}=\frac {3{c}^{2}-{b}^{2}-2{c}^{2}b\pm\sqrt {4{c}^{4}{b}^{2}-
12{c}^{4}b-4{c}^{2}{b}^{3}+9{c}^{4}+2{c}^{2}{b}^{2}+{b}^{4}}}{2cb}.
\end{align*}

Based on linear analysis, the type of  equilibria for   \eqref{2d} (\(a=0\)) and for $b,c\in [-4,4]\times [-2,2]$ in different regions defined by \eqref{2d} are provided in Figures \ref{1fig} and \ref{CASEA}. \\
 
\begin{enumerate}
\item We  
prove the existence of a Pitchfork bifurcation  for \(b=1\) using the normal form truncated at order 3  \cite[Section 20.1E]{wiggins}. Let \(b=1+\mu\),  where \(\mu \in \RR\) is an additional bifurcation parameter.  For \(\mu=0\),  the Jacobian matrix \eqref{JE2_DZb} associated with \(E_1\) admits eigenvalues: \(0\) and \(\frac{c^2-1}{c}\), $c\neq 0,1$. Define the diffeomorphism $T: \RR^2\to \RR^2$ as
\begin{equation*}
\mathbf{x}\equiv(x,y) \mapsto  \left(\frac{c^2 y-x}{c^2-1}, -\frac{y-x}{c^2-1}\right).
\end{equation*}
In the new coordinates $\mathbf{x}=(x,y),$ where $\0=(0,0)$, let $$\bar{f}_{(0, b,c)}(\mathbf{x})\equiv f_{(0, b,c)}(\mathbf{x})-Df_{(0, b,c)}(\0)\mathbf{x}.$$ Therefore, we have 
 \begin{align}\label{r1}
\dot{\mathbf{x}}= 
&T^{-1} Df_{(0, b,c)}(\0) T\mathbf{x}+T^{-1}\bar{f}_{(0, b,c)}T\mathbf{x}\\\nonumber
=&\begin{pmatrix}
1 & c^2\\
1 & 1
\end{pmatrix}
\begin{pmatrix}
 c & c\\ 
-\frac{1}{c} &-\frac{1}{c}
\end{pmatrix}
\begin{pmatrix}
\frac{1}{1-c^2} &  \frac{c^2}{c^2-1}\\ 
\frac{1}{c^2-1} & \frac{1}{1-c^2} 
\end{pmatrix}
+ \\
&\begin{pmatrix}
-c y\mu-\frac{1}{3} c x^3 \\\nonumber
\frac{c^2-1}{c}x +\frac{c^2-\mu-1}{c} y-\frac{c}{3} x^3
\end{pmatrix}\Big|_{x\rightarrow \frac{c^2 y-x}{c^2-1}, y\rightarrow \frac{x-y}{c^2-1}}
\\\nonumber
\end{align}
from where we conclude that (in the new coordinates above):
\begin{equation*}
 \left\{ 
\begin{array}{l}
\dot{x}=-\frac{c}{3(c^2-1)^3} 
\Big(3c^2x^2y-x^3-3c^4xy^2+c^6y^3+3(c^2-1)^2\mu x -3 (c^2-1)^2 \mu y \Big),\\ \\ 
\dot{y}=
-\frac{1}{c(c^2-1)^3} \Big(3c^4x^2y-c^2x^3-3c^6xy^2+c^8y^3+3 (c^2-1)^2 \mu x-3 (c^2-1)^2 \mu y-3 (c^2-1)^4 y\Big).
\end{array}
\right.
\end{equation*}

We aim to find a center manifold in the form of 
$$h(x,\mu) := {\alpha_0}\mu^2+\alpha_1 x^2+\alpha_2 x\mu + \dots.$$ Since the center manifold must satisfy the equation (cf. \cite[Equation 3.2.7]{wiggins})
$$\dot{x} \frac{\partial h}{\partial x}-\dot{y}\big|_{(x,h(x,\mu))}=0,$$ we obtain \(\alpha_0=\alpha_1=0,\) and \(\alpha_2=\frac{1}{\left(c^2-1\right)^2}.\)
 Then, the dynamics of   \eqref{r1} restricted to the center manifold 
 $y=\frac{\mu x}{(c^2-1)^2}$ is given by:
\begin{equation}
\label{pitchfork1.eq}
\dot{x}=-\frac{c}{c^2-1}\mu x+\frac{c}{3\left(c^2-1\right)^3} x^3
\end{equation}
meaning that there is a Pitchfork bifurcation  at $\mu=0$. If $c>1$, then the bifurcation is supercritical with respect to $\mu$ and also with respect to $b$ in the original equation, as depicted in Figure \ref{pitchfork1}. The analysis for the case $c>1$ is different from the case $0<c<1$, although the argument runs along the same lines. \\
 
\begin{figure}[t!]
\begin{center}
 \includegraphics[height=4.9cm]{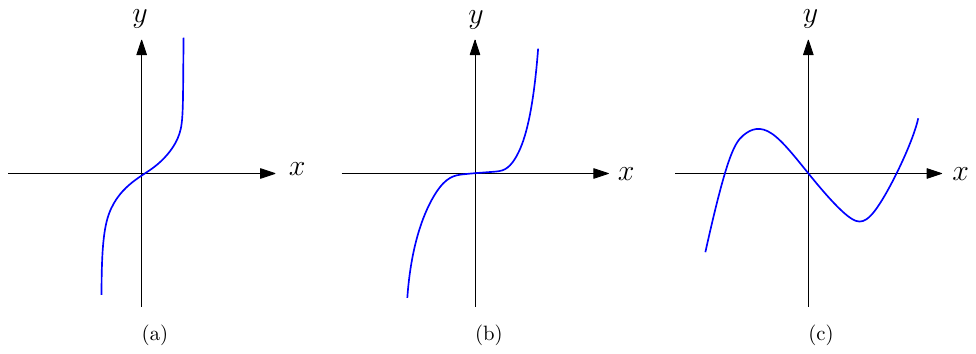}
  \end{center}
\caption{ \label{pitchfork1}    \small   Illustration of the Pitchfork bifurcation of   \eqref{pitchfork1.eq} for  $c>1$ and close to $\mu=0 (\Leftrightarrow b=1)$.  (a) $\mu<0$, (b) $\mu=0$ and (c) $\mu>0$.  }
\end{figure}

\item   We should prove that the eigenvalues of the Jacobian matrix of $f_{(0,b,c)}$ at $E_1$ change from real to complex (non-real).  The proof of this item follows by observing that the eigenvalues associated with the Jacobian matrix \eqref{JE2_DZb} at \(E_1\)  are complex (non-real) for 
\(b > -c^2 + 2c\) or \(b < -c^2 - 2c,\) and real otherwise. The sign of their real part does not change.\\

\item Let $b=c^2$ and $|c| <1$. The Jacobian matrix of $f_{(0,b,c)}$ at $E_1$ has  a pair of purely imaginary eigenvalues \(i \sqrt{1-c^2}\) with \(-1<c<1.\) Let \(b=c^2+\mu\) where \(\mu\in \RR\) is a bifurcation parameter. The eigenvalues associated with the Jacobian matrix at $E_1$ for $\mu\neq  0$   
\begin{equation*}
\begin{pmatrix}
c & c \\ 
-\frac{1}{c} & -c-\frac{1}{c}\mu
\end{pmatrix},
\end{equation*}
are given by
\begin{eqnarray}
\lambda(\mu)=-\frac{1}{2c}\mu\pm \frac{1}{2c}\sqrt{(\mu+2c^2+2c) (\mu+2c^2-2c)}.
\end{eqnarray} 
 
For $\mu$ close to zero and $|c|<1$, there exists a pair of complex conjugate eigenvalues. Furthermore, since 
\begin{eqnarray*}
\frac{d \, {\rm Re}{\lambda(\mu)}}{d\mu}=-\frac{1}{2c}\neq 0,
\end{eqnarray*} 
this means that there is a non-degenerate Hopf bifurcation at $\mu=0$ (\cite[pp. 299--300]{Meiss}). By using the linear transformation 
\begin{align*}
\begin{pmatrix}
u \\ 
v
\end{pmatrix} =
\begin{pmatrix}
\dfrac{1}{2} & \dfrac{c}{2\sqrt{1-c^2}} \\ 
0 & -\dfrac{1}{2c\sqrt{1-c^2}}
\end{pmatrix} 
\begin{pmatrix}
x \\ 
y
\end{pmatrix},
\end{align*}
the system \eqref{2d} turns into its Jordan form 
\begin{equation*}
 \left\{ 
\begin{array}{l}
f^1: \, \dot{u}= -\frac{1}{12}c{u}^{3}-\frac{1}{4}\frac {{c}^{2}}{\sqrt {1-{c}^{2}}}v{u}^{2}
-\frac{1}{4}{\frac {{c}^{3}}{1-{c}^{2}}}{v}^{2}u
-\frac{1}{12}\frac {{c}^{4}}{\left( 1-{c}^{2} \right) ^{3/2}}{v}^{3}-\sqrt {1-{c}^{2}}v,\\ \\ 
f^2: \, \dot{v}= \sqrt {1-{c}^{2}} u.\end{array}
\right.
\end{equation*}

The first Lyapunov coefficient  depends on $c $ and is given by:
\begin{align*}
\ell(c)= \dfrac{1}{16\omega} \left(R_1+ \omega R_2\right),
\end{align*}
where
\begin{align*}
 \omega:=&\sqrt {1-{c}^{2}} \neq 0 \\
R_1:=&f^1_{uv} \left(f^1_{uu} + f^1_{vv} \right)
-f^2_{uv} \left(f^2_{uu} + f^2_{vv} \right)
 -f^1_{uu} f^2_{uu}+f^1_{vv} f^2_{vv}\quad \text{and}\\\nonumber
R_2:=&f^1_{uuu} + f^1_{uvv} + f^2_{uuv} + f^2_{vvv}.
\end{align*}

\bigbreak

Using Maple, we get $$\ell(c)=\frac{c}{32 (c^2-1)}<0$$ and we may conclude that the limit cycle is stable for  \(0<c<1\) and the bifurcation is 
supercritical
 -- \cite[Remark 1, pp. 383]{wiggins}.
In this case, $E_1$ is an asymptotically stable equilibrium for $\mu > 0$ ($\Leftrightarrow b>c^2$)  and  unstable  for $\mu < 0$ ($\Leftrightarrow b<c^2$), with an asymptotically stable periodic solution for $\mu < 0$ ($\Leftrightarrow b<c^2$). \\

\item The proof is similar as item (2) by taking into account that\\
\begin{eqnarray*}
&4{c}^{4}{b}^{2}-
12{c}^{4}b-4{c}^{2}{b}^{3}+9{c}^{4}+2{c}^{2}{b}^{2}+{b}^{4}&\\\nonumber
&= \left(b^2+2cb-2c^2b+3c^2\right)\left(b^2-2cb-2c^2b+3c^2\right).&
\end{eqnarray*}
\\

\item Using linear analysis, one knows that Hopf bifurcation exists only if 
\(\mbox{Tr}(E_2) = \mbox{Tr}(E_3)=0\) and \(\mbox{Det}(E_2), \mbox{Det}(E_3)>0.\) Indeed, we have: \\
\begin{eqnarray*}
&& \mbox{Tr}(E_2) = \mbox{Tr}(E_3)=0  \mbox{ and } b>1\\
& \Leftrightarrow& 3{c}^{2}-{b}^{2}-2{c}^{2}b =0  \mbox{ and } b>1\\
&\Leftrightarrow&  b=c\left(-c\pm\sqrt{c^2+3}\right) \mbox{ and } b>1. \\
\end{eqnarray*}

\begin{figure}[t]
\begin{center}
 {\includegraphics[width=.68\columnwidth,height=.65\columnwidth]{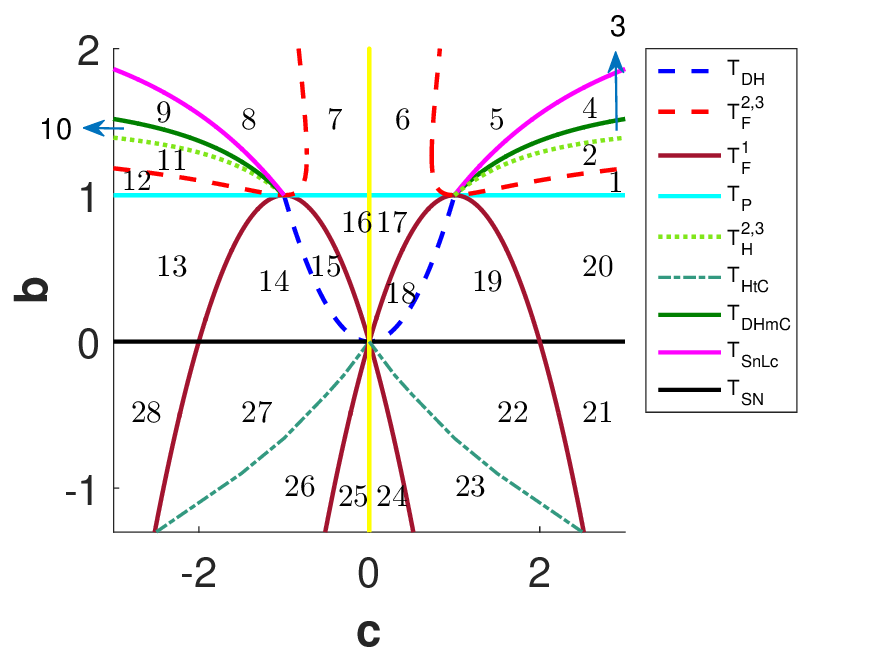}}
\end{center}
\caption{\label{CASEA} The bifurcation curves associated to  \eqref{2d} for \textbf{Case A}.    The local portraits corresponding to Regions \textbf{1--28} have been plotted in Table \ref{QA} of Appendix.} 
\end{figure}

The inequality \(b>1\) serves to ensure that the equilibria $E_2$ and $E_3$ exist. \\ \\

Let \(b=c\left(-c\pm\sqrt{c^2+3}\right).\) The Jacobian matrix at \(E_{2,3}\) has  a pair of purely imaginary eigenvalues 
\(\pm i\sqrt{ 2 c\sqrt{c^2+3}-2-2 c^2}\) with \(c>1.\) Let \(b=-c^2\pm c\sqrt{c^2+3}+\mu\) where $\mu\in \RR$ denotes the bifurcation parameter. The eigenvalues, \(\lambda(\mu),\) associated with the Jacobian matrix at $E_{2,3}$
\begin{eqnarray*}
\begin{pmatrix}
-\frac{c \left( 2{c}^{2}+2c\sqrt {{c}^{2}+3}\pm2\mu\mp3 \right) }{\mp{c}^{2}+c\sqrt {{c}^{2}-3}\pm \mu}&c\\
-\frac{1}{c} &{\frac {{c}^{2}\mp c\sqrt {{c}^{2}+3}-\mu}{c}}
\end{pmatrix},
\end{eqnarray*}
are  given by (computations performed in Maple):
 \begin{eqnarray*}
&\lambda(\mu)=\dfrac{-2 c \sqrt{c^2+3}\mu-\mu^2\pm\sqrt{\sigma}}{2c (\mu-c^2+c\sqrt{c^2+3})}&
\end{eqnarray*}
where
\begin{small}
\begin{eqnarray*}
&\sigma=\mu^4+4c(\sqrt{c^2+3}-2c)\mu^3-4c^2(6c\sqrt{c^2+3}-7c^2-5)\mu^2&\\\nonumber
&+8c^3(2\sqrt{c^2+3} (3 c^2+1)-6c^3-11c)\mu
-8c^4(c\sqrt{c^2+3}\left( 4c^2+5\right)-11c^2-4c^4-3).&
\end{eqnarray*}
\end{small}
 
  There is a pair of complex conjugate eigenvalues when \(\sigma<0\). In addition, we have: 
$$\frac{d\,  {\rm Re}{\lambda(\mu)}}{d\mu}=-\frac{\sqrt{c^2+3}}{c(\sqrt{c^2+3}\mp c)}.$$
The first Lyapunov coefficient  $\ell$ associated to   $f_{(0,b,c)}$ at $E_{2,3}$ depends on $c$ and is given explicitly by (computations in Maple):
\begin{eqnarray*}
&\ell(c)={\frac {c \left(  \left( 4c^{2}+3 \right)  \left( \sqrt {c^{2}+3}-c \right) -3c \right)  
\left( \sqrt {c^{2}+3}c \left( 21+44c^{2}+16c^{4} \right) -c^{2} 
\left( 69+68c^{2}+16c^{4} \right) -9 \right) }
{32 \left( \sqrt {c^{2}+3}c-c^{2}-1\right) ^{2} \left( \sqrt {c^{2}+3}-c \right) ^{7}}},&
\end{eqnarray*}
 
 Since $\displaystyle \lim_{c \to 1^+} \ell(c) = \lim_{c \to +\infty} \ell (c) = +\infty$ and
$ \displaystyle \frac{\partial \ell (c )}{\partial c} = 0 \quad \text{at} \quad c^\star \approx 1.315$, then $\ell$ remains positive in the neighbourhood of  $c^\star$. Therefore, the corresponding limit cycles for \(c > 1\) are unstable and the bifurcation is subcritical.\\

\item This item follows from \cite[pp. 201]{GlobalFN} and \cite[Formula (3.2.29)]{fitzhughbook}.
\end{enumerate}
\end{proof}

\begin{rem}
Numerically, we may observe that the non-degenerate Hopf bifurcation  of items (3) and (5) of Theorem \ref{th:main1} are, respectively, 
 supercritical and subcritical.
The bifurcation $T_H^1$ yields a stable periodic solution, say $C^s$, and $T_H^{2,3}$ generates two unstable and $\kappa$--symmetric periodic solutions -- confirm on Table \ref{QA} of Appendix.\\ 
\end{rem}

\subsection{Case B (\(b=0\))}
 
Model \eqref{2d} for \(b=0\) has only one equilibrium given explicitly by: 
\begin{eqnarray*}
E_1= \,\left(a,-a+\frac{1}{3} a^3\right),
\end{eqnarray*}
 and the linear part of \(f_{(a,0, c)}\) at $E_1$ has eigenvalues
 \begin{eqnarray}\label{eigb0}
\frac{1}{2} c(1-a^2)\pm \frac{1}{2}\sqrt{(c a^2-c+2) (ca^2-c-2)}.
\end{eqnarray}

Then we may conclude that: \\

\begin{thm}
\label{Th3.8}
With respect to the vector field $f_{(a,0,c)}$ of \eqref{2d},  the global representations of the transition/bifurcation curves in the 
parameter space
\((a,c)\in \RR\times \RR\backslash \{0\}\) are as follows:\\
\begin{itemize}
\item[I.]
The equilibrium \(E_1\) undergoes a Belyakov transition along the hyperbola: 
\begin{align*}
T^{b0}_F=\left\{(a,c): c=\pm \frac{2}{a^2-1}\right\}. \\
\end{align*}
 
 \begin{figure}
\begin{center}
  {\includegraphics[width=.65\columnwidth,height=.56\columnwidth]{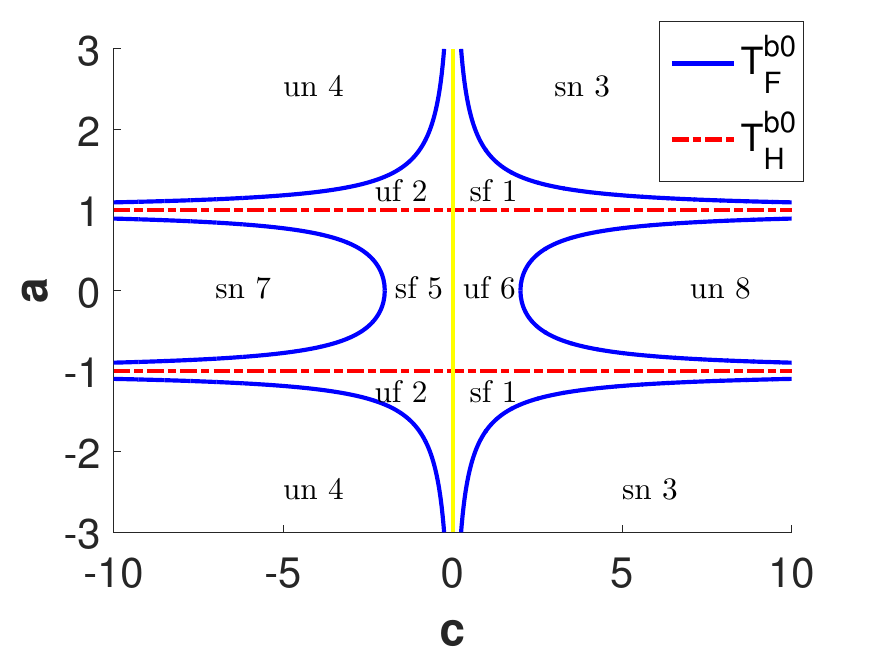}}
 \qquad
 \caption{\label{CASEC}  {  The classification of the equilibrium point $E_1$ and the corresponding bifurcation curves for   \eqref{2d} in \textbf{Case B}. In this figure, \textbf{un, sn, uf} and \textbf{sf} refer to unstable (stable) node and unstable (stable) focus, respectively.  The local phase portraits associated with Regions \textbf{1--8} have been plotted in Table~\ref{QB}. }}
\end{center}
\end{figure}

\item[II.]
The equilibrium \(E_1\) undergoes  a Hopf bifurcation along the lines: 
 \begin{align*}
T^{b0}_{H}=\left\{(a,c): a=\pm 1\right\}.
\end{align*}
The Hopf bifurcation is supercritical(resp. subcritical) if $c>0$ (resp. $c<0$).\\
\end{itemize}
\end{thm}
Schematic bifurcation lines for \(b=0\) have been plotted  in Figure \ref{CASEC}. The local phase portraits associated with Regions \textbf{1--8} have been plotted in Table~\ref{QB}
\bigbreak
\begin{proof}
 We are going to concentrate the proof in the case $b=0$ and $c>0$.   Again, the analysis of the case $b=0$ and $c<0$   has a similar treatment (reversing time).\\
\begin{itemize}

\item[I.] 
  The first claim is a direct consequence of    \eqref{eigb0}.\\
\item[II.]
 
 The Jacobian matrix at \(E_{1}\) has  a pair of purely imaginary eigenvalues 
\(\pm i.\) Let \(a=\pm1+\mu\) where \(\mu \in \RR\) refers to the bifurcation parameter. The eigenvalues, \(\lambda(\mu),\) corresponding to the Jacobian matrix \eqref{JE2_DZb} at \(E_{1}\) are (computations performed with Maple):
 \begin{eqnarray*}
&\lambda(\mu)=\mp c\mu-\frac{1}{2} c {\mu}^{2}\pm \frac{1}{2}\sqrt {4{c}^{2}{\mu}^{2}\pm 4{c}^{2}{\mu}^{3}+{c}^{2}{\mu}^{4}-4}. &
\end{eqnarray*}
For $\mu$ close to zero, we have $4{c}^{2}{\mu}^{2}+4{c}^{2}{\mu}^{3}+{c}^{2}{\mu}^{4}-4<0$, which means that $\lambda(\mu)$ is complex non-real. Furthermore 
$$\frac{d \, {\rm Re}{\lambda(\mu)}}{d\mu}\big|_{\mu=0}=\mp c\neq 0.$$ 
The first Lyapunov coefficient depends on $c$ and is given by $$\ell(c)=-\frac{1}{32}c$$ which has opposite sign to that of $c\neq 0$.  
Hence the bifurcation is subcritical for $c<0$ and supercritical for $c>0$.
 
 For the case $a=-1$,   $E_1$ is an asymptotically stable equilibrium for $\mu < 0$ and  unstable focus  for $\mu > 0$ with an asymptotically stable periodic solution emerging for  $\mu > 0$.  For the case $a=1$,   $E_1$ is an asymptotically stable equilibrium for $\mu > 0$ and  unstable focus  for $\mu < 0$ with an asymptotically stable periodic solution emerging $\mu > 0$.  
This finishes the proof. \end{itemize}
\end{proof}

\begin{thm}\label{thm11}
With respect to the vector field $f_{(a,0,c)}$ of \eqref{2d}, if $|a|<1$ then there exists at most one limit cycle. If it exists, then it is unstable.
\end{thm}

For proving Theorem \ref{thm11} we shall use the next theorem:

\begin{thm}[Cherkas and Zhilevich
\cite{cherkas}, adapted] \label{main12}
Consider the system of differential equations
\begin{equation} \label{2d_intro}
  \left\{ 
\begin{array}{l}
\dot{x}=   -\phi(y) - F(x), \\ \\
 \dot{y}=  \displaystyle g(x)
\end{array}
\right.
\end{equation}
 where
\[
F(x) = \int_0^x f(x)\,dx.
\]

Suppose that for every $x \in (a,b)$, with
$-\infty \le a < 0 < b \le +\infty$, and every
$-\infty < y < +\infty$, the following conditions are satisfied:\\

\begin{enumerate}
    \item[(i)] $xg(x) > 0$ for $x \ne 0$,
    $y\phi(y) > 0$ for $y \ne 0$;

    \item[(ii)] $f(x)$, $g(x)$, $\phi(y)$ are continuously differentiable,
    $\phi(y)$ is monotonically increasing,
    $f(0) < 0$ (or $f(0) > 0$);

    \item[(iii)] there exist constants $\alpha$, $\beta$ such that
    \(
    f_1(x) = f(x) + g(x)\left[\alpha + \beta F(x)\right]
    \)
    has simple zeros $x_1 < 0$ and $x_2 > 0$, and
    $f_1(x) \le 0$ (or $f_1(x) \ge 0$)
    in the interval $(x_1,x_2)$;

    \item[(iv)] the function $f_1(x)/g(x)$ is nondecreasing
    (or nonincreasing) outside the interval $[x_1,x_2]$;

    \item[(v)] every closed orbit encloses the interval $[x_1,x_2]$
    on the $x$-axis.\\
\end{enumerate}

Then,  \eqref{main12} can have at most one limit cycle.
The limit cycle is stable (or unstable), if it exists.

\end{thm}

\begin{proof}
The differential system \eqref{2d} when $b=0$ has the unique equilibrium point
\(
\left(a,\frac{a(a^2-3)}{3}\right).
\)
Doing the translation
\(
(x,y) = \left(X+a,\; Y+\frac{a(a^2-3)}{3}\right)
\)
we send this equilibrium to the origin of coordinates of:
\begin{eqnarray}
\label{main13}
&& \dot{X} = cY + c(1-a^2)X - acX^2 - \frac{c}{3}X^3,\\
\nonumber&& \dot{Y} = -\frac{X}{c}.
 \end{eqnarray}

Now we rescale the time $t$ as 
\(
ds = -c\,dt.
\)
In the new time $s$ the differential system \eqref{main13} becomes
\begin{eqnarray}
\label{main14}
&&X' = \frac{dX}{ds}
= -Y - (1-a^2)X + aX^2 + \frac{1}{3}X^3,\\
\nonumber&&Y' = \frac{X}{c^2}.
 \end{eqnarray}

The differential system \ref{main14} is written in the form \eqref{main12} with

\begin{eqnarray*}
\label{main14}
&&\phi(Y)=Y, \\
&&F(X)=(1-a^2)X-aX^2-\frac{1}{3}X^3,\\
\nonumber&&g(X)=\frac{X}{c^2}.\\
 \end{eqnarray*}

Since
\(
Xg(X)=\frac{X^2}{c^2}$ and $Y\phi(Y)=Y^2, $ then 
the condition (i) of Theorem \ref{main12} holds. \\

Following the notation of Theorem \ref{main12} we have
\[
f(X)=1-a^2-2aX-X^2,
\]
and consequently
\(
f(0)=1-a^2,
\)
distinct from $0$ if $a \ne \pm 1$.
Hence if $a \ne \pm 1$ the condition (ii) of Theorem \ref{main12} is satisfied.

Taking $\alpha = \beta = 0$ we have that
\(
f_1(X)=f(X),
\)
and the function $f_1(X)$ has the two simple zeros
\(
X_1=-1-a,\)
and
\(X_2=1-a.\)
Therefore
\(
f_1(X) \ge 0
\)
in the interval $(X_1,X_2)$.
So if $|a|<1$ the condition (iii) of Theorem \ref{main12} holds. Since the function
$$
\frac{f_1(X)}{g(X)}
= -c^2\frac{(X+a-1)(X+a+1)}{X}
$$
is nonincreasing outside the interval $[X_1,X_2]$,
the condition (iv) of Theorem \ref{thm11} is satisfied. \\

The divergence of the differential system \eqref{main14} is
\(
(X-1+a)X + (1+a).
\)
Since the divergence does not change sign inside the strip
\(
\{(X,Y): X_1 < X < X_2\},
\)
it follows that every closed orbit encloses the interval
$[x_1,x_2]$ on the $X$-axis.
Hence condition (v) of Theorem \ref{main12} holds.\\

In summary, the five conditions of Theorem \ref{main12} are satisfied by the differential system \eqref{2d} 
with $b=0$ and $|a|<1$.
The result follows straightforwardly (note that $f(0)>0$). \\
\end{proof}

\subsection{Case C (\(a > 0,\) \(0<b<1\))}


If \((x^\star, y^\star)\) is an equilibrium of (\ref{2d}) under the conditions  \(a > 0,\) \(0<b<1,\) then:
\begin{eqnarray*}
x^\star - a +by^\star=0 \Leftrightarrow  y^\star = \frac{a - x^\star}{b}.
\end{eqnarray*}
 
By substituting \(y^\star = \frac{a - x^\star}{b}\) into the first equation of \eqref{2d}, then \(x^\star\) is a root of \(H(x)=0,\) where
\begin{eqnarray*}
H(x) = -\frac{1}{3} c x^3 + c \left(1 - \frac{1}{b}\right) x + \frac{a c}{b}.
\end{eqnarray*}
 
The discriminant associated with the polynomial \(H(x)\) is given by 
\begin{eqnarray*}
\Delta = \frac{4}{3} c^4 \left(1 - \frac{1}{b}\right)^3 - 3 \frac{c^4 a^2}{b^2}.
\end{eqnarray*}
For \(0 < b < 1\), the discriminant satisfies \(\Delta < 0\).  According to Cardano's formula, this implies that the polynomial $H(x)$ has exactly one real root, denoted by \(x^\star.\)
Then the equilibrium for system (\ref{2d}) is explicitly given by 
\begin{eqnarray*}
\left(x^\star,   \frac{a - x^\star}{b}\right)=(x^\star, y^\star)=:E_1.
\end{eqnarray*}
The characteristic equation of the Jacobian matrix evaluated at \((x^*, y^*)\) is 
\begin{eqnarray*}
\lambda^2 + \frac{c^2 {x^*}^2 - c^2 + b}{c} \lambda + b {x^*}^2 - b + 1 = 0.
\end{eqnarray*}
We may then conclude that:\\

\begin{prop}
\label{Prop3.10}
The following assertions hold for the vector field $f_{(a,b,c)}$ of \eqref{2d},  under the conditions \(a > 0,\) and \(0<b<1\): 
\begin{enumerate}
\item 
The equilibrium point $E_1$ is stable if \(\frac{c^2 {x^*}^2 - c^2 + b}{c} > 0.\)
Conversely, it is unstable if \(\frac{c^2 {x^*}^2 - c^2 + b}{c} < 0.\)
\item 
A non-degenerate Hopf bifurcation occurs at  \(c^2 = \frac{b}{1 - {x^*}^2}.\)
The numerical simulations  indicate the resulting limit cycle is attracting when \(c > \sqrt{\frac{b}{1 - {x^*}^2}}\) and repelling for \(c < -\sqrt{\frac{b}{1 - {x^*}^2}}\).\\
\end{enumerate}
  \end{prop}
Note that $x^\star \neq \pm 1$ (if  $x^\star = \pm 1$ then $b=0$, which is a contradiction). 
 
As before, 
 we use the notation $\mbox{Det}(E_{1})$ and $\mbox{Tr}(E_{1})$ for the determinant and trace
operators of \(Df_{(a,b, c)},\)  evaluated at $E_1$.

\bigbreak
\begin{proof}
 
Under the condition \(0<b<1,\) one knows that $\det(E_1)>0.$\\
\begin{enumerate}
\item The result follows by noticing that $\mbox{Tr}(E_{1})=- \frac{c^2 {x^*}^2 - c^2 + b}{c} .\) If $\mbox{Tr}(E_{1})<0$, then the  equilibrium $E_1$ is stable.  Conversely, it is unstable if $\mbox{Tr}(E_{1})>0$.\\ 
\item 
 The Jacobian matrix at \(E_{1}\) has  a pair of purely imaginary eigenvalues 
\(\pm i \sqrt{b {x^*}^2 - b + 1}.\) Let \(b=c^2 \left(1 - {x^*}^2\right)+\mu\) where \(\mu\in \RR\) refers to the bifurcation parameter. The eigenvalues, \(\lambda(\mu),\) corresponding to the Jacobian matrix at \(E_{1},\) 

\begin{eqnarray*}
\begin{pmatrix}
c(1-x^*) & c\\
 -\frac{1}{c} & c{x^*}^2-c-\frac{\mu}{c}
\end{pmatrix},
\end{eqnarray*}
are 
 \begin{eqnarray*}
&\lambda(\mu)=\frac{-\mu\pm \sqrt{\left(2c^2 {x^*}^2-2 c^2-\mu-2 c\right) \left(2c^2 {x^*}^2-2 c^2-\mu+2 c\right) }}{2c}.&
\end{eqnarray*}
If 
  $$\min_{c\geq 0}{c\left(\pm1-c+c {x^*}^2\right)}<\mu<\max_{c\geq 0}{ c\left(\pm1-c+c {x^*}^2\right)},$$  
 then the expression under the square root is negative. Hence there exists a pair of complex non real conjugate eigenvalues. 
 The result follows by observing that (computations performed with Maple): \[\frac{d {\rm Re}\, {\lambda(\mu)}}{d\mu}=-\frac{1}{2c}\neq 0.\]
\end{enumerate}
\end{proof}
 
 \begin{rem}
 The statement of Proposition \ref{Prop3.10} says nothing about the criticality of the Hopf Bifurcation. Analytically, it is very difficult to compute the first 
 Lyapunov coefficient of \eqref{2d} at an equilibrium whose explicit expression we do not know. However, numerics   suggest that the Hopf bifurcation is supercritical for $c=   \sqrt{\frac{b}{1 - {x^*}^2}}$ and subcritical for  $c= -  \sqrt{\frac{b}{1 - {x^*}^2}}$. \\
 \end{rem}
 
 \section{Analysis of the singular case \(c \to+\infty\) in Cases A and B}
\label{s:canards}
The 
aim
of this section is the study of the dynamics of \eqref{2d} when $|c| \to +\infty$. 
 Many  results concerning the asymptotic behaviour of the solution of the FN model as $|c|\to +\infty$ are available in the literature -- see \cite{fitzhughbook, 2D, connections}.  Our concern is to match these results with the local bifurcation study performed in Section \ref{sec3}.   We are particularly interested in \textbf{Cases A} and \textbf{B} whose results have been described in Theorems \ref{th:main1} and \ref{Th3.8}.
 By changing the time scale of \eqref{2d}, 
 taking into account
 \(\tau= t/c\) and \(\varepsilon = 1/c^2,\)  
 we get an equivalent form
for
 \eqref{2d} as follows (written with respect to the \emph{slow time} $\tau$):
\begin{equation} 
\label{2d_epsilon}
\left\{ 
\begin{array}{l}
\varepsilon \dot{x}= \displaystyle \left[y-\left(\frac{x^3}{3}-x\right)\right]=: f(x,y, \varepsilon) \\ \\
 \dot{y}=  \displaystyle - (x+by-a)=: g(x,y, \varepsilon)
\end{array}
\right.
\end{equation}
where \(a, b\in \RR\) and \(c\in \RR\backslash\{0\}.\)
When \(0<\varepsilon \ll 1,\) equivalently  \(|c|\gg  1,\) system \eqref{2d_epsilon} 
is termed
a \emph{fast-slow system} (see \cite{GLR2025, Kuehn} and references therein).

A fast-slow system consists of differential equations where the rates of change of some variables are much greater than those of others.
 This leads to a two time-scale system.
 The general approach to this type of system starts by grouping the variables in two disjoint sets: fast variable $x$ and slow variable $y$. 
 
By decomposing the dynamics into slow and fast phases, each can be studied independently, allowing a clearer understanding of its underlying behaviour.
 The parameter \(\varepsilon\)  introduces this separation in system \eqref{2d_epsilon}.
 
\subsection{The critical manifold and bifurcations}
\begin{figure}[ht]
\begin{center}
\subfigure[]
 {\includegraphics[width=.55\columnwidth,height=.46\columnwidth]{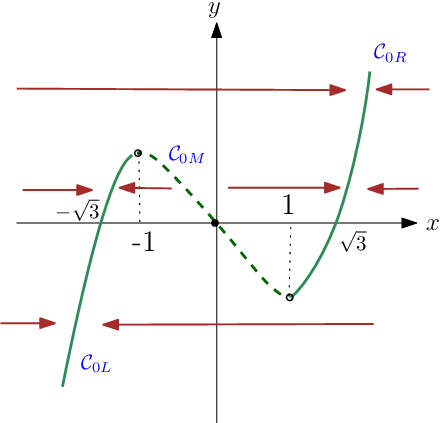}}
 \qquad
 \subfigure[]
 {\includegraphics[width=.55\columnwidth,height=.46\columnwidth]{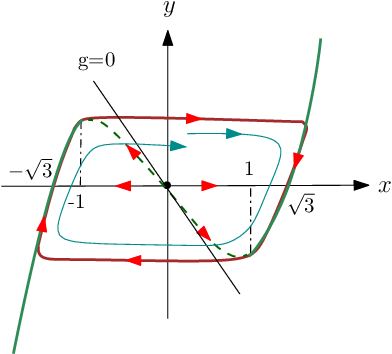}}
  \end{center}
 
\caption{   \label{slow-fast1} \small (a)  Illustration of the critical manifold associated to  system \eqref{2d_epsilon}.  (b)   Illustration of the dynamics of system \eqref{2d_epsilon} for the case \(a = 0,\) when there is a unique equilibrium.}
\end{figure}

  We first analyse the singular case $\varepsilon=0.$
We follow the analysis performed in \cite{GLR2025, KS2001a, KS2001b}. The critical manifold associated with \eqref{2d_epsilon} is the set 
\begin{eqnarray}
\label{C0_def}
\mathcal{C}_0=\{(x,y) \in \RR^2: f(x, y, 0)=0\}= \left\{(x,y) \in \RR^2: y= \frac{x^3}{3}-x\right\}.
\end{eqnarray}
 \bigbreak
 \begin{defn}
 \label{nh_def}
  With respect to \eqref{2d_epsilon}  with critical manifold $ \mathcal{C}_0$:
        \begin{enumerate} 
            \item The subset $S\subset \mathcal{C}_0$ is called  \emph{normally hyperbolic} provided that for all $e=(x_0,y_0)\in{S}$, we get:
                                \begin{equation*}
              \frac{\partial f}{\partial x}\left(e,0\right) \neq 0
            \end{equation*}
\medskip
        
 \item A normally hyperbolic set $S$ is called  \emph{attracting} (resp. \emph{repelling}) 
 when,  for every $e\in{S}$, all eigenvalues of $  \frac{\partial f}{\partial x}\left(e,0\right)$ possess negative (resp. positive) real part, respectively. 
 
        \end{enumerate}
 \end{defn}
  \bigbreak
  \begin{defn} \label{def_fold}
         With respect to \eqref{2d_epsilon} with critical manifold $\mathcal{C}_0$ (see \eqref{C0_def}), the point $e\in \mathcal{C}_0$ 
          is called
         a  \emph{fold point} if:
        \begin{equation*}
             \frac{\partial f}{\partial x}\left(e,0\right)=0\,,\quad  \frac{\partial^2 f}{\partial x^2}\left(e,0\right)\neq0, \quad  \frac{\partial f}{\partial y}\left(e,0\right)\neq0.
        \end{equation*}
        If $f\left(e,0\right)\neq0$, the fold point is termed  \emph{regular}.
        
    \end{defn}
  \bigbreak
Coming back to system \eqref{2d_epsilon}, 
we get $f\left(x,y,0\right)= x-\frac{x^3}{3}+y$ and, for $e=(x_0,y_0)\in \mathcal{C}_0$, one has: 
\begin{align*}
\frac{\partial f}{\partial x}\left(e,0\right)= 1-x_0^2\neq0\qquad\text{ if }\,\, x_0\neq\pm 1.
\end{align*}
        Therefore, if $b=3/2$,  the equilibria $E_2=\left(-1,2/3\right)$ and $E_3=(1,-2/3)$   are fold points (according to the Definition  \ref{def_fold})  since:
        $$
        \frac{\partial f}{\partial x}(E_{2,3},0)=0, \quad \frac{\partial^2 f}{\partial^2 x}(E_{2,3},0)\neq 0\quad \text{and}\quad \frac{\partial f}{\partial y}(E_{2,3},0)\neq 0
        $$
                   and $S_0=\mathcal{C}_0\setminus\left\{E_2, E_3\right\}$ is normally hyperbolic according to Definition \ref{nh_def}.
        We may split $S_0$ into three disjoint open subsets:
\begin{equation*}
\begin{aligned}
                \mathcal{C}_{0\text{L}}&= \mathcal{C}_0 \cap {\textstyle \left\{\left(x,y\right)\in\RR^2: x<-1 \right\}}, \\
                \mathcal{C}_{0\text{M}}&= \mathcal{C}_0 \cap {\textstyle \left\{\left(x,y\right)\in\RR^2:-1<x<1\right\}}\quad \text{and} \\
                \mathcal{C}_{0\text{R}}&= \mathcal{C}_0 \cap {\textstyle \left\{\left(x,y\right)\in\RR^2: x>1 \right\}}. \\
            \end{aligned}
        \end{equation*}
        Therefore, we get 
        \begin{equation}
           \label{def: S0}
    S_0=\mathcal{C}_{0\text{L}}\cup\mathcal{C}_{0\text{M}}\cup\mathcal{C}_{0\text{R}},    
        \end{equation}
where  $\mathcal{C}_{0\text{L}}$ and $\mathcal{C}_{0\text{R}}$ are \emph{attracting}  and $\mathcal{C}_{0\text{M}}$ is \emph{repelling}. 
As a result, trajectories of \eqref{2d_epsilon} in the fast flow 
go 
towards either to $\mathcal{C}_{0\text{L}}$ or $\mathcal{C}_{0\text{R}}$ and backwards to $\mathcal{C}_{0\text{M}}$ -- see Figure \ref{slow-fast1}.\\

 Yet in the singular case ($\varepsilon=0$), depending on the values of $a,b\in \RR$ in system  \eqref{2d_epsilon}, 
 one, two, or three equilibria might arise as a consequence of the cubic curve's intersection, \(f (x, y, 0) = 0,\)  and 
the line defined by $x-a+by=0$ ($\Leftrightarrow \dot{y} = 0).$ 
If the system possesses three equilibria, only one of the following cases   arise
(cf.  \cite{GLR2025}): \\ 

\begin{description}
\item[Case 1]  
 
each of the three equilibria lies in  $\mathcal{C}_{0\text{M}}$ or else 
\item[Case 2] 
 each of the equilibria belongs to a different region $\mathcal{C}_{0\text{L}}$, $\mathcal{C}_{0\text{M}}$, and $\mathcal{C}_{0\text{R}}$.\\
\end{description}

When \(\varepsilon = 0,\) the independence of the time scales simplifies the system's analysis.
 Fenichel's Theorem \cite{Fenichel} says that, near the normally hyperbolic part of  $\mathcal{C}_{0}$, 
 for
\(0<\varepsilon\ll1,\) 
the perturbed system exhibits dynamics ``analogous'' to the singular system  (in the sense of differentiable conjugacy),
with a deviation of the order \(\mathcal{O}\left(\varepsilon\right).\) 
 Decreasing \(\varepsilon\) makes 
the trajectories of the system increasingly similar to those of the singular system ($\varepsilon=0$).
 According to Fenichel's Theorem, it is ensured that
when
\(0<\varepsilon\ll1,\) 
  there is a set \(S_\varepsilon\), $\mathcal{O}\left(\varepsilon\right)$--close to $S_0$ (in the \(\mathbb{C}^1\)--topology),  which is locally hyperbolic and attracting (if $S_0$ is normally hyperbolic and attracting).
 \begin{defn}
The set $S_\varepsilon$ is known as the \emph{slow manifold} for \eqref{2d_epsilon}.
  \end{defn}

When the equilibrium of the system is a fold point,
a subcritical (with respect to $b$) 
 Hopf bifurcation arises, leading to the emergence of a periodic solution.
In this section, the \emph{canard phenomenon} described in \cite{KS2001a, KS2001b} 
links
 the local periodic solution 
 generated by
 the Hopf bifurcation to the global limit cycle, 
 like
 we proceed to explain (the canard point is correlated with the passage of a reduced problem equilibrium over the fold point under parameter variation -- see Definition \ref{def:regular_fold}).
\begin{defn}
           A trajectory segment of a fast-slow system  \eqref{2d_epsilon} is a   \emph{{canard}} if it stays within a $\mathcal{O}\left(\varepsilon\right)$--distance to a repelling branch of a slow manifold for a time that is $\mathcal{O}\left(1\right)$ on the slow time scale $\tau = t\varepsilon$.           
\end{defn}
          
\emph{Canards} are related to equilibria of    \eqref{2d_epsilon} that arise at fold points of the slow manifold. The next definitions deal with the general equation parametrised by $\lambda\in\RR$:
\begin{equation}
\label{general1} 
    \left\{ 
\begin{array}{l}
\varepsilon \dot{x}=   f(x,y, \lambda, \varepsilon) \\  \\
 \,\dot{y}=  g(x,y, \lambda, \varepsilon)
\end{array}
\right.
\end{equation}

   We now define a \emph{singular fold point},
   a different concept from     that of a fold point; 
     it is necessary to
 with an equilibrium point of  \eqref{2d_epsilon}:
    \begin{defn}\label{def:fold_point}
    With respect to \eqref{general1}, let $e=\left(x,y\right)$ 
       denote
     a fold point. 
     We refer to $e$ as a \emph{singular fold} provided:
        \begin{equation}
            \label{dobrasingular}
            \begin{aligned}
                &f\left(e,\lambda,0\right)=0, &&\frac{\partial f}{\partial x}\left(e,\lambda,0\right)=0,\\
                &\frac{\partial ^2 f}{\partial x^2}\left(e,\lambda,0\right)\neq0, && \frac{\partial f}{\partial y}\left(e,\lambda,0\right)\neq0, \quad g\left(e,\lambda,0\right)=0.
            \end{aligned}
        \end{equation}     
    \end{defn}

    \begin{defn}\label{def:regular_fold}
        A singular fold point $e$ (also known as a canard point) 
                 is called
         \emph{regular} if:
        \begin{equation}
            \frac{\partial g}{\partial x}\left(e,\lambda,0\right)\neq0 \quad \text{ and } \quad \frac{\partial g}{\partial \lambda}\left(e,\lambda,0\right)\neq0.
        \end{equation}
    \end{defn}

    The dynamics near a regular singular fold point 
     is explored  in Theorem \ref{th:canard}.
     The occurrence of \emph{canards} is linked to a Hopf bifurcation, known in these systems as a \emph{singular Hopf bifurcation}.     
    
 \begin{thm}[\cite{KS2001a, KS2001b}, adapted]
   \label{th:canard}
        Consider a  fast-slow system of the form \eqref{general1} where $\left(x,y\right)=\left(0,0\right)$ is 
             canard point
        for $\lambda=0$, and 
               \begin{equation}\label{eq:nformHopfsing}
            \begin{aligned}
                &\dot{x}=-y\,l_1\left(x,y,\lambda,\varepsilon\right)+x^2\,l_2\left(x,y,\lambda,\varepsilon\right) + \varepsilon\,l_3\left(x,y,\lambda,\varepsilon\right)\\
                &\dot{y}=\varepsilon \left(\pm x\,l_4\left(x,y,\lambda,\varepsilon\right)-\lambda \,l_5\left(x,y,\lambda,\varepsilon\right)+y\,l_6\left(x,y,\lambda,\varepsilon\right)\right)
            \end{aligned}
        \end{equation}
        where
        $$
      l_3\left(x,y,\lambda,\varepsilon\right)= \mathcal{O} \left(x,y,\lambda,\varepsilon\right) \quad \text{and}\quad l_j\left(x,y,\lambda,\varepsilon\right)= 1+\mathcal{O} \left(x,y,\lambda,\varepsilon\right), 
        $$
        for $j\in\{1,2,4,5,6\}$.
 
          Suppose that, when $\varepsilon=0$, a slow trajectory connects attracting and repelling regions of the critical manifold $\mathcal{C}_0$.
        Therefore, there is $\varepsilon_0>0$ and $\lambda_0>0$ so that for $0<\varepsilon<\varepsilon_0$ and $|\lambda|<\lambda_0$, the system possesses an equilibrium point $p\in\RR^2$ around 
        \((0,0)\)
                 where \(p\) approaches \((0,0)\) as $\left(\lambda,\varepsilon\right)$ tends to \((0,0)\).
Moreover,  there 
 is
a smooth map $\lambda_c:\left[0,\varepsilon_0\right]\to\RR$ 
   that 
assigns
every
 $\varepsilon\in \left(0,\varepsilon_0\right]$ to  a real value $\lambda$,
leading to a family of \emph{canards}, 
defined asymptotically as
        \begin{equation*}
            \lambda_c\left(\sqrt{\varepsilon}\right)=-\varepsilon\left(A + B\right)+\mathcal{O}\left(\varepsilon^{3/2}\right).  
        \end{equation*}
        and 
in addition, there is a continuous  map $\lambda_ {H}:\left[0,\varepsilon_0\right]\to\RR$
  that associates 
 every
  $\varepsilon\in \left[0,\varepsilon_0\right]$ to a value $\lambda$, associated to a Hopf bifurcations  
 defined asymptotically as
        \begin{equation*}
            \lambda_ {H}\left(\sqrt{\varepsilon}\right)=
            -\varepsilon B\varepsilon+\mathcal{O}\left(\varepsilon^{3/2}\right),    
        \end{equation*}
        where
        \begin{equation*}
            A=-\dfrac{\dfrac{\partial l_1}{\partial x}-3\dfrac{\partial l_2}{\partial x}+2\dfrac{\partial l_4}{\partial x}-2l_6}{8},\qquad B=\dfrac{\dfrac{\partial l_3}{\partial x}+l_6}{2},    
        \end{equation*}
        %
 and the functions $l_i$  and their partial derivatives, for $1\leq i\leq 6$,
   are calculated at $\left(x,y,\lambda,\varepsilon\right)=\left(0,0,0,0\right)$.  
  The latter bifurcation 
          is non-degenerate when $A \neq 0$, it occurs for  
   positive values of $\lambda-\lambda_H$  if $A < 0$, and  for  
negative  values of $\lambda-\lambda_H$ if $A > 0$.
       \end{thm}

\subsection{Case A \((a=0)\)}    
 
  The curves \(T_F^{2,3}\), \(T_H^{2,3}\), $DH$ and SNL, introduced in Theorem~\ref{th:main1} and Remark~\ref{rem:SNL}, can be  parameterised by \(c \in [1, +\infty)\). From now on, we restrict our attention to the portions of these curves lying within the region defined by \(b, c \geq 1\) -- see Figure \ref{summary1}.
  We are going to denote by $T_{F, +} ^{2,3} \backslash\{(1,1)\}$ the connected component of $T_F^{2,3} \backslash\{(1,1)\}$ with highest values of $b$ for the same $c$, and $T_{F, -} ^{2,3} \backslash\{(1,1)\}$ the other.
   It is straightforward to verify that, for sufficiently large values of \(c\) in \([1, +\infty)\), the following holds:
\begin{eqnarray}
T_{F, -}^{2,3} (c)\leq   T_H^{2,3}(c) \leq DH(c) \leq  SNL(c)  
.\\ \nonumber
\end{eqnarray}
All the curves pass through $(b,c)=(1,1)$. Denote by \(A_1 \leq A_2<A_3<A_4\) the limits  when $c\to +\infty$ of the previous curves (when it exists), respectively. More specifically, we have:
 
 \begin{figure}
\begin{center}
\includegraphics[height=10.5cm]{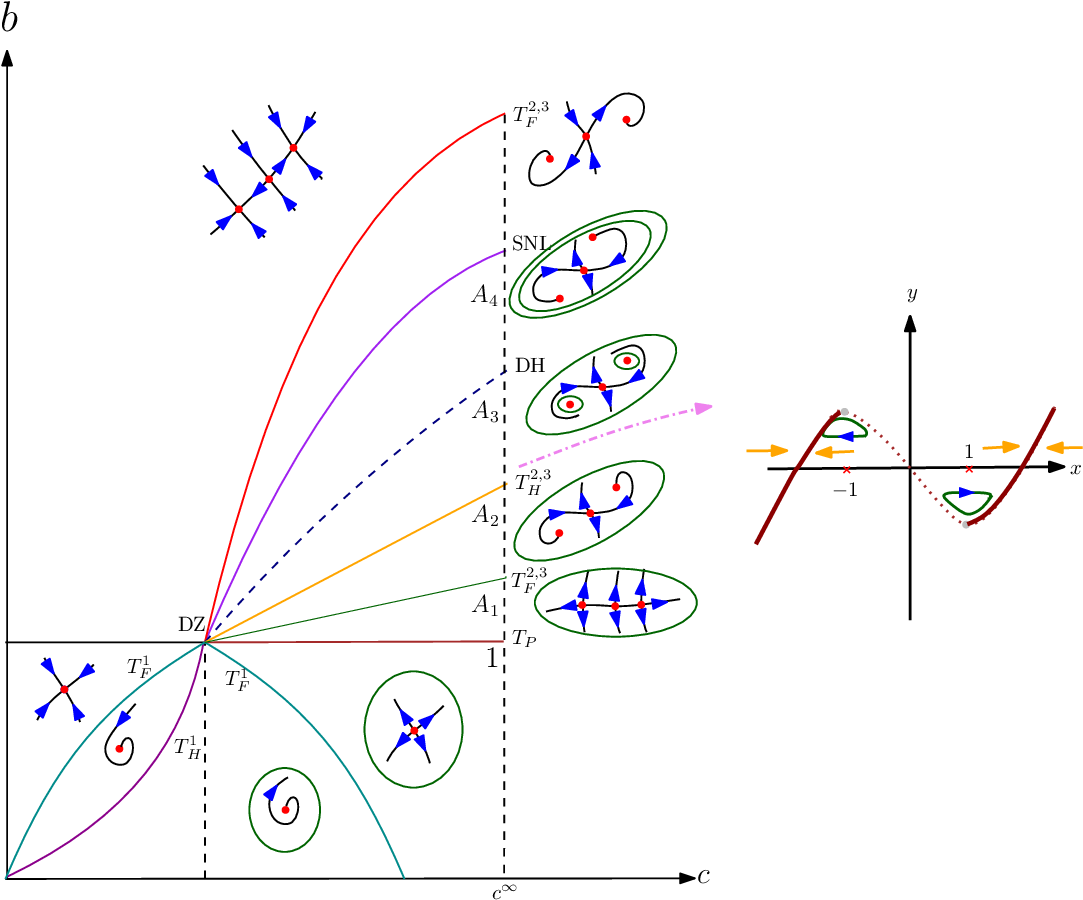}
 \caption{ Partial scheme of the codimension 1 bifurcation curves emerging from the point $(b,c)=(1,1)$ associated with $f_{(0,b,c)}$ of (\ref{2d}), for $b,c> 0$. The symbol $c^\infty$ represents ``very large'' values of $c$. The notation of the curves correspond to those of Theorem \ref{th:main1}, for \textbf{Case A}.  The line SNL is dashed because it does not make part of the DZ bifurcation with $\mathbb{Z}_2(\kappa)$--symmetry. Near the curve $b=T_H^{2,3}$, the emerging periodic solutions give rise to  canards as a consequence of Proposition \ref{canardsA}. }
 \label{summary1}
\end{center}
\end{figure}
 
 \begin{lem} As illustrated in Figure \ref{summary1}, the following equalities hold:
\begin{enumerate}
\item \(A_1= \displaystyle \lim_{c\to +\infty}   T_{F, -}^{2,3}(c)=3/2 \)
\item \(A_2= \displaystyle \lim_{c\to +\infty}   T_H^{2,3}(c)= 3/2\)
\item \(A_3=\displaystyle \lim_{c\to +\infty}   DH(c)= {17/10}\)
\item \(A_4= \displaystyle \lim_{c\to +\infty}  SNL(c)= +\infty \)
\end{enumerate}
In particular, we have $1<A_1\leq A_2<A_3<A_4.$
\end{lem}
We omit the proof since it follows from straightforward computations using the expressions of Theorem \ref{th:main1} and Remark \ref{rem:SNL}. A partial representation of these curves has been plotted in Figure \ref{summary1}.

\bigbreak

 \begin{prop}
 \label{canardsA}
With respect to \eqref{general1}, 
a smooth function $b_c:\left[0,\varepsilon_0\right]\to\RR$  exists
that 
assigns
every
$\varepsilon\in \left(0,\varepsilon_0\right]$ to a value $b$ 
leading to
 a family of \emph{canards} in the system, asymptotically defined by: 
        \begin{equation*}
            b_c\left(\sqrt{\varepsilon}\right)=\frac{3}{2}+\frac{5}{4} \varepsilon+\mathcal{O}\left(\varepsilon^{3/2}\right).   
        \end{equation*}     
 \end{prop}
 
 \begin{proof}
 The difference $b-3/2$ of system \eqref{general1} is played by $\lambda$ of Theorem  \ref{th:canard}. 
The equilibrium \(E_3\) is a singular fold point and  
a Hopf bifurcation occurs at $b=3/2$.
Using
Theorem \ref{th:canard}, we perform a change of variables and parameter as
\begin{equation*}
\left\{ 
\begin{array}{l}
\bar{x}=x-1\\ \\
\bar{y}= y+2/3 \\ \\
\lambda=b-3/2 \\
\end{array}
\right.
\end{equation*}
to obtain the equivalent system ($\tau=t/\varepsilon$)
\begin{equation*}
\left\{ 
\begin{array}{l}
\dot{\bar{x}}=\bar{y} -\bar{x}^2-\frac{1}{3}\bar{x}^3\\ \\
\dot{\bar{y}}=\varepsilon\left(-\bar{x}  +\frac{2}{3}\lambda   -\lambda \bar{y}  - \frac{3}{2}  \bar{y}\right).
\end{array}
\right.
\end{equation*}
 
Using the notation introduced in Theorem \ref{th:canard}, we obtain
\begin{eqnarray*}
l_1=-1\qquad
l_2= -1-\bar{x}/3\qquad
l_3=0\qquad
l_4=1\qquad
 l_6=-\frac{3}{2}
\end{eqnarray*}
hence \(A=- 1 /2<0\) and $B= -3/4<0$ therefore the Hopf bifurcation occurs for values $\lambda>\lambda_H (\Leftrightarrow  b>\frac{3}{2}+\lambda_H )$ . Since for $b<T_H^{2,3}(c)$ the eigenvalues of  \(Df_{(0,b, c)},\)  evaluated at \(E_i,\) \(i=1,2,\) are negative, then the bifurcating periodic solution is unstable and canards occur when 
$$
b_c(\sqrt{\varepsilon})=\frac{3}{2}+\frac{5}{4} \varepsilon + \mathcal{O}(\varepsilon^{3/2}).
$$       
\end{proof} 
   
\subsection{Digestive remark}  
 We describe the dynamics of  (\ref{2d}) in the case \(a = 0\), focusing on the first quadrant of the bifurcation diagram in the \((b, c)\)--plane. 
 We suggest the reader follows the description by using Figure \ref{summary1}.
 
For \(b \in (0, 1)\) and \(c >  1\), there exists a unique equilibrium point \(E_1\), which is unstable. This source is enclosed by a stable periodic solution denoted by \(C^s\).
At \(b = 1\), one of the eigenvalues of $Df_{(0,b,c)}$ at  \(E_1\) becomes zero, leading to the creation of two additional  unstable nodes \(E_2\) and \(E_3\). In other words,  \(E_1\) undergoes a supercritical (with respect to $b$) Pitchfork bifurcation.

The saddles $E_2$ and $E_3$ subsequently become unstable foci and undergo a supercritical Hopf bifurcation, 
leading to
two unstable periodic solutions (symmetric under $\kappa$).
Around \(b = DH(c)\), the invariant manifolds of $E_1$ evolve into a double homoclinic orbit, which is unstable, previously identified in~\cite{GlobalFN}.

\begin{figure}[ht]
\begin{center}
 \includegraphics[height=1.4cm]{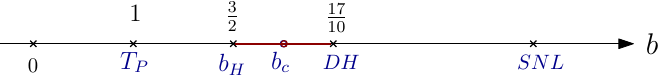}
  \end{center}
\caption{ \label{Axis} Sketch of the singular bifurcations for $b\geq 1$ and $c<\infty$, in \textbf{Case A}.  
As $\varepsilon$ decreases, the interval $[b_c, b_H]$ becomes narrower.
When $\varepsilon =0 \, (\Leftrightarrow c=+\infty$), we have $b_c\equiv b_H$. }
\end{figure}

 As \(b\) increases smoothly, the double homoclinic loop breaks apart, resulting in an unstable periodic orbit \(C^u\), which surrounds both the stable and unstable manifolds of \(E_1\). 
Near \(b = SNL(c)\), the periodic orbits \(C^s\) and \(C^u\) collapse, and two stable equilibria remain, dominating the basin of attraction of any compact subset of \(\mathbb{R}^2\). For $c \gg1$, there is a small portion of the bifurcation parameter diagram $(b,c)$ of Figure \ref{summary1}, where the periodic orbits emerging from $T_H^{2,3}$ correspond to canards.   They rapidly collapse into the double homoclinic cycle to $E_1$. 
Parameters $b_c$ that yield \emph{canards} are only $\mathcal{O}\left(\varepsilon\right)$ away from the corresponding Hopf bifurcation values $b_{H}$
(see Figure \ref{Axis}):
    \begin{equation*}
        b_{c}-b_H=\frac{1}{2}\,\varepsilon+\mathcal{O}\left(\varepsilon^{3/2}\right).    
    \end{equation*}
 
  \subsection{Case B \((b=0)\)}    
  \begin{prop}
  \label{canardsB}
With respect to \eqref{general1}, 
a continuous mapping $a_c:\left[0,\varepsilon_0\right]\to\RR$ exists
that assigns to each $\varepsilon\in \left[0,\varepsilon_0\right]$ a corresponding value $a$ 
giving rise to
a \emph{canard},
defined asymptotically by
        \begin{equation*}
            a_c\left( \sqrt{\varepsilon}\right)=-\frac{\varepsilon}{8}+\mathcal{O}\left(\varepsilon^{3/2}\right).   
        \end{equation*}
       
 \end{prop}

\begin{proof}
The difference $a-1$ of system \eqref{general1} is played by $\lambda$ in Theorem  \ref{th:canard}. 
We are interested in only one of the two fold points of (\ref{2d}) located at $(1, -2/3)$.  After shifting the fold point to the origin, reversing time $t \to - t$ and setting $\lambda= a-1$, we get
\begin{equation*}
\left\{ 
\begin{array}{l}
\dot{x}=- y + x^2 ( 1 + x/3 ),\\ \\
\dot{y}=\varepsilon (x-\lambda).
\end{array}
\right.
\end{equation*}  

In this standard form, the relevant parameters $\ell_i$ defined in Theorem \ref{th:canard} are $$\ell_1=-1, \quad \ell_2=1+x/3, \quad \ell_3=0, \quad \ell_4=1, \quad \ell_5=1,\quad  \ell_6=0, \quad A=1.$$
Hence, a maximal canard exists on a curve $\lambda_c$ in $(\varepsilon, \lambda)$--space defined by:
$$
a_c(\sqrt{\varepsilon})= -\frac{\varepsilon}{8} + \mathcal{O}(\varepsilon^{3/2}).
$$
\end{proof}
  \subsection{Digestive remark}
 
 The dynamics of system \eqref{2d} for \(b = 0\) is illustrated in Figures \ref{CASEC} and Table \ref{QB}.
Our analysis focuses on the first and fourth quadrants of the bifurcation diagram in the \((a, c)\)--plane and $c>0$.
For \(a < -1\), the saddle \(E_1\) becomes a stable node.
Near \(a = -1\), a Belyakov transition occurs, followed by a supercritical Hopf bifurcation that gives rise to a stable periodic orbit that involves into a canard in an exponentially small portion of the bifurcation plane.  This periodic solution disappears through another Hopf bifurcation at \(a = 1\).

\section{The phase portraits ``near infinity'' } \label{sec4}
  Studying the dynamics ``near infinity'' provides insight into the system's asymptotic behaviour and allows the identification of  trajectories that leave bounded sets.
Using the techniques of  Subsections \ref{sec:vertical} and \ref{ss:compactification(prel)}, we compactify system  \eqref{2d} and   describe the  dynamics of  \eqref{2d} at infinity -- this is what we call the  dynamics ``near infinity".   We distinguish between four cases
\begin{align}
b\geq 0, \,c>0, \qquad b\geq 0, \,c<0, \qquad b< 0, \,c>0, \qquad b<0, \,c>0.
\end{align}

 The main theorem of this section is the following:

\begin{thm}\label{S50_th}
The dynamics of system \eqref{2d} ``near infinity'' is topologically equivalent to one of the phase portraits 
depicted in Figure \ref{ID}.
\end{thm}

\subsection*{Insight into the reasoning}
Before proceeding with the formal proof of the result, we give a geometrical interpretation of what is going to be done:  
we investigate the behaviour of trajectories of the planar vector field \eqref{2d} as their position tends toward infinity. To do this, we apply the Poincar\'e compactification of Subsection  \ref{ss:compactification(prel)}, which transforms points at infinity into a bounded region   represented by a disk containing a finite number of equilibria.
   We then study the qualitative dynamics of the ``compactified'' vector field    paying particular attention to the equilibria located there, the directions along which trajectories diverge, and the conditions under which the flow at infinity becomes attracting or repelling. Some equilibria are singularities. In order to study their dynamics we need to desingularize them using the techniques of Subsection \ref{sec:vertical}.\\
\begin{figure}[t]
\begin{center}
\subfigure[{\tiny \(b\geq 0\) and \(c>0\)}]
 {\includegraphics[width=.22\columnwidth,height=.22\columnwidth]{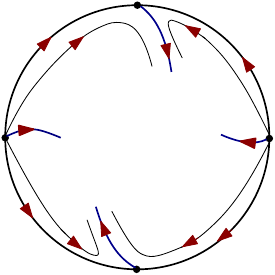}}
 \quad
 \subfigure[{\tiny  \(b\geq 0\) and \(c<0\)}]
 {\includegraphics[width=.22\columnwidth,height=.22\columnwidth]{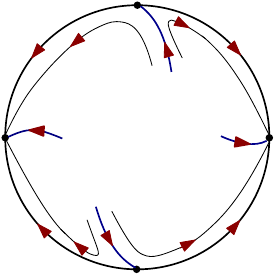}}
 \quad
 \subfigure[{\tiny  \(b<0\) and \(c<0\)}]
 {\includegraphics[width=.22\columnwidth,height=.22\columnwidth]{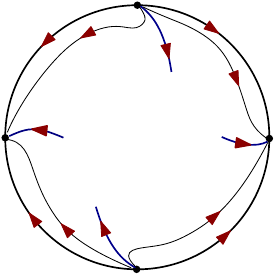}}
    \quad
 \subfigure[{\tiny  \(b<0\) and \(c>0\)}]
  {\includegraphics[width=.22\columnwidth,height=.22\columnwidth]{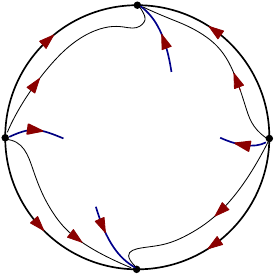}}
  \end{center}
  \caption{\label{ID}  Dynamics near infinity on the Poincar\'e disc of system \eqref{2d}.}   
\end{figure}

\begin{proof}
Using \eqref{comp1} with $d=3$, system \eqref{2d} in the local chart  \(U_1\)
may be written as:
\begin{align}\label{XXX}
\dot{u}=&-cu^2v^2-\frac{b}{c}uv^2-cuv^2+\frac{a}{c}v^3-\frac{1}{c}v^2+\frac{c}{3}u,\\\nonumber
\dot{v}=&-cuv^3-cv^3+\frac{c}{3}v,
\end{align}
where
$$
 u= \frac{x_2}{x_1}\quad\text{and}\quad  v=\frac{1}{x_1}.
$$
 \bigbreak

The system \eqref{XXX} has a unique equilibrium at the origin
on \(v=0\) 
(in any chart, points of \(\mathbb{S}^1\) are characterized by \(v= 0\)).
The eigenvalues of the Jacobian matrix associated to \eqref{XXX} computed at the origin,   are \(\dfrac{c}{3}, \dfrac{c}{3}.\) So, the origin is an unstable (resp. stable) node for \(c>0\) (resp. \(c<0\)). \\
Analogously, using \eqref{comp2} and $d=3$, the system \eqref{2d} in chart \(U_2\)
takes the form:
\begin{align}\label{k8}
\dot{u}=&\frac{1}{c}u^2v^2-\frac{a}{c}uv^3-\frac{c}{3}u^3+\frac{b}{c}uv^2+cuv^2+cv^2,\\\nonumber
\dot{v}=&\frac{1}{c}uv^3-\frac{a}{c}v^4+\frac{b}{c}v^3.
\end{align}
The origin is an identically zero singular point for system \eqref{k8}. In order to desingularize this equilibrium, we are going to employ the quasi-homogeneous blow-up
\begin{align}\label{tr1}
(u, v) = (\bar{u} \rho^{2}, \bar{v} \rho^{3}),
\end{align}
where \((\alpha,\beta)=(2,3)\) has been obtain by using the Newton diagram   discussed in Subsection \ref{sec:vertical}.  The Newton polygon associated to \eqref{k8} is illustrated in Figure \ref{Newton1}.

 \begin{figure}
\begin{center}
\includegraphics[height=5.5cm]{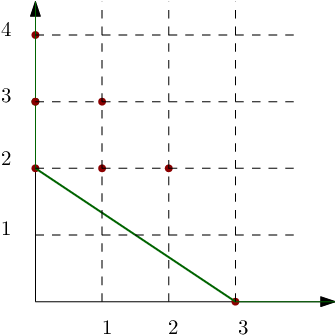}
\end{center}
\caption{\small The points/segments associated to \eqref{k8} are $(2,2)$, $(1,3)$, $(3,0)$, $(1,2)$, $(0,2)$, $(0,4)$ and $(0,3)$. The Newton polygon is defined by the line $3y+2x\geq 6$. Therefore  $(\alpha,\beta)=(2,3)$.}
 \label{Newton1}
\end{figure}

Applying the transformation \eqref{tr1} where \(\bar{u} = 1\) and \(\bar{u}=-1,\) 
the \(x\) directional blow-up, converts the system \eqref{k8} into the vector fields defined by: 
\begin{align*}
X_1^{+}:
\begin{cases}
\frac{d{\rho}}{d\eta} =\frac{\rho \left(-3a \rho^5{\bar{v}}^{3}+ 3{\rho}^{4}{\bar{v}}^{2}+3{c}^{2}{\rho}^{2}{\bar{v}}^{2}
+3b{\rho}^{2}{\bar{v}}^{2}+3{c}^{2}{\bar{v}}^{2}-{c}^{2} \right) }{6c}
,\\\\\nonumber
\frac{d\bar{v}}{d\eta}  =-\frac {\bar{v} \left(-a \rho^5{\bar{v}}^{3} +{\rho}^{4}{\bar{v}}^{2}+3{c}^{2}{\rho}^{2}{\bar{v}}^{2}
+b{\rho}^{2}{\bar{v}}^{2}+3{c}^{2}{\bar{v}}^{2}-{c}^{2} \right) }{2c},
\end{cases}
\end{align*}
and 
\begin{align*}
X_1^{-}:
\begin{cases}
\frac{d{\rho}}{d\eta} =\frac{\rho \left(-3a \rho^5{\bar{v}}^{3} -3{\rho}^{4}{\bar{v}}^{2}+3{c}^{2}{\rho}^{2}{\bar{v}}^{2}
+3b{\rho}^{2}{\bar{v}}^{2}-3{c}^{2}{\bar{v}}^{2}-{c}^{2} \right) }{6c}
, \\\\\nonumber
\frac{d\bar{v}}{d\eta} =-\frac{\bar{v} \left(-a \rho^5{\bar{v}}^{3}+ 5{\rho}^{4}{\bar{v}}^{2}+3{c}^{2}{\rho}^{2}{\bar{v}}^{2}
+b{\rho}^{2}{\bar{v}}^{2}+3{c}^{2}{\bar{v}}^{2}-{c}^{2} \right) }{2c},
\end{cases}
\end{align*}
respectively, where $d\eta = \rho^{4} d\tau$ (the old time variable is $\tau$ and the new one is $\eta$; when we calculate $X^1$, we can cancel $\rho^4$ by using this time rescaling). \\

We now need to analyze the equilibria at \(\rho=0\)  of \(X_1^{\pm}.\)
Regarding the vector field \(X_1^+,\) and  \(X_1^-,\) one concludes that: \\

\begin{enumerate}
\item[\textbf{Claim 1:}]
The equilibria of system \(X_1^{+}\) on \(\rho = 0\) are directly determined by the roots of 
$$\bar{v}  \left(3 \bar{v}^{2} -1\right)=0.$$ 
There are three equilibria at \(\rho = 0\) denoted by$$D_1: (0,0)\qquad \text{and} \qquad f^{\pm}: \left(0, \pm\dfrac{\sqrt{3}}{3}\right).$$ 
The following statements are valid:\\
\begin{itemize}
\item[1.] The equilibrium \(D_1\) always is a saddle point for \(X_1^{+}\).\\
\item[2.] The points \(f^{\pm}\) are semi-hyperbolic and  the dynamics of   \(X_1^{\pm}\) restricted to the center manifold 
 $\bar{v}=\mp\frac{\sqrt{3} (b+3c^2) \rho^2}{18 c^2}$  is given by \(\dot{\rho}=\frac{b}{9c}{\rho^3}.$
 
 \begin{proof}
 A center manifold in the form of 
$h(\rho) := \alpha_1 \rho^2+ h.o.t$ must satisfy the equation 
\[\dot{\rho} \frac{\partial h}{\partial \rho}-\dot{\bar{v}}\big|_{(\rho,h(\rho))}=0.\] 
After shifting \(f^{\pm}\) to the origin, we obtain \(\alpha_1=\mp\frac{\sqrt{3} (b+3c^2)}{18 c^2}.\)
 Then, the dynamics of   \(X_1^{\pm}\) restricted to the center manifold 
 $\bar{v}=\mp\frac{\sqrt{3} (b+3c^2) \rho^2}{18 c^2}$ 
 is given by \(\dot{\rho}=\frac{b}{9c}{\rho^3}\) 
 (more details in \cite[Theorem 2.19]{llibrebook}).
 \end{proof}
 
\item[3.] The equilibria  \(f^{\pm}\)  are: \\
\begin{itemize}
\item[1a.] saddle points for$
\{(a,b,c) \big| b\geq 0, \, c<0\} \cup \{(a,b,c) \big| b\geq 0, \, c>0\}$ -- see Figures \ref{cca} and \ref{ccb}.\\
\item[2a.] an unstable node for
$
\{(a,b,c) \big| b<0, \, c<0\}$ -- see Figure \ref{cce}.\\
\item[3a.] a stable node for
 $\{(a,b,c) \big| b<0, \, c>0\}.
$ -- 
see Figure \ref{ccd}.\\
\end{itemize}

\end{itemize}
 
\item[\textbf{Claim 2:}]
The equilibria of system \(X_1^{-}\) on \(\rho = 0\) are determined by the roots of 
$$\bar{v}  \left(3 \bar{v}^{2} -1\right)=0.$$ There are three equilibria at \(\rho = 0\) denoted by \(D_2: (0,0)\) and \(g^{\pm}: \left(0, \pm\frac{\sqrt{3}}{3}\right),\) and the results of Claim 1 are valid in this case (changing $D_1$ and $f^\pm$ by $D_2$ and $g_\pm$, respectively).\\

Using the transformation \eqref{tr1} with \(\bar{v}=1\) and \(\bar{v}=-1,\) the \(y\) directional blow-up  transforms system \eqref{k8} into
\begin{align*}
Y_1^{+}:
\begin{cases}
\frac{d{\rho}}{d\eta} =\mathcal{O}(|\rho, \bar{u}|), \\\\
\frac{d\bar{u}}{d\eta} =c,
\end{cases}
\quad\quad
Y_1^{-}:
\begin{cases}
\frac{d{\rho}}{d\eta} =\mathcal{O}(|\rho, \bar{u}|), \\\\
\frac{d\bar{u}}{d\eta} =c,
\end{cases}
\end{align*}
respectively, where \(d\eta = \rho^{4} d\tau.\) 
 We now need to analyze the point (0, 0) of \(Y_1^{\pm}.\) Indeed, \\ 
 
\item[\textbf{Claim 3:}]For the vector fields \(Y_1^{\pm},\) point (0, 0) is a regular point (not an equilibrium).\\
\end{enumerate}

\begin{figure}[t!]
\begin{center}
\subfigure[{\tiny \(b\geq 0\) and \(c>0\)}]
 {\includegraphics[width=.22\columnwidth,height=.22\columnwidth]{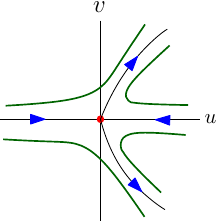}}
 \quad
 \subfigure[{\tiny  \(b\geq 0\) and \(c<0\)}]
 {\includegraphics[width=.22\columnwidth,height=.22\columnwidth]{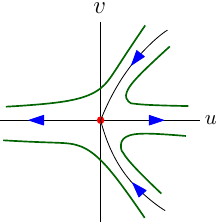}}
 \quad
 \subfigure[{\tiny  \(b<0\) and \(c<0\)}]
 {\includegraphics[width=.22\columnwidth,height=.22\columnwidth]{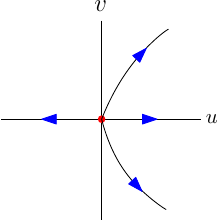}}
    \quad
 \subfigure[{\tiny  \(b<0\) and \(c>0\)}]
  {\includegraphics[width=.22\columnwidth,height=.22\columnwidth]{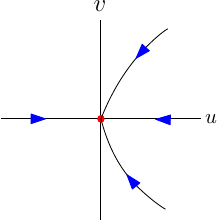}}
  \end{center}
  \caption{\label{ccc}  Local phase portraits of system \eqref{k8}  (the representation of system \eqref{2d} in the chart \(U_2\)) at \((0,0)\).} 
\end{figure}

  From the \textbf{Claims 1--3}, we know the dynamics of the vector fields \(X_1^{\pm}\) and \(Y_1^{\pm}\). For the vector fields \(Y_1^{\pm},\) \((0, 0)\) is a regular point. Thus, Figure \ref{cc} illustrates the qualitative properties of the vector fields \(X_1^{\pm}\) and \(Y_1^{\pm}.\) By shrinking the circles  in these figures to the origin (as the radii of the circles decrease to zero and the circles collapse at 
  \((0,0)\)
  ), we 
  achieve to
   the phase portraits shown in Figure \ref{ccc}. So, as shown in Figure \ref{ID}, the equilibrium of chart \(U_2\) is a saddle with four hyperbolic sectors for \(b\geq0.\) Moreover,  it becomes an unstable (resp. stable) node when \(b,c<0\) and \(b<0, c>0,\) respectively.
 
\begin{figure}[t!]
\begin{center}
\subfigure[{\tiny \(b\geq 0\) and \(c>0\) } \label{cca}]
 {\includegraphics[width=.22\columnwidth,height=.22\columnwidth]{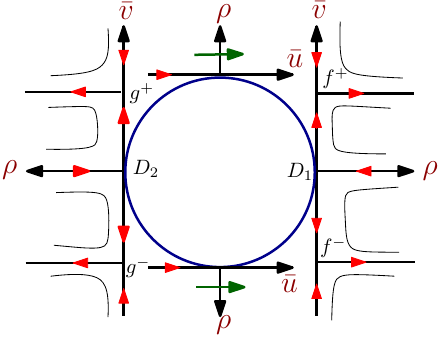}}
 \quad
 \subfigure[{\tiny \(b\geq 0\) and \(c<0\)} \label{ccb}]
 {\includegraphics[width=.22\columnwidth,height=.22\columnwidth]{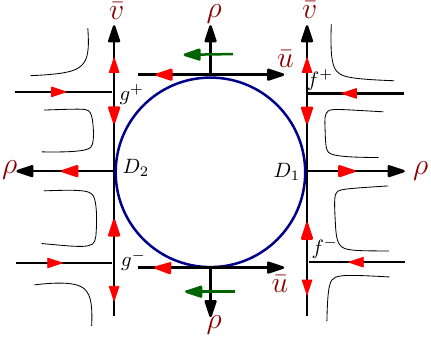}}
 \quad
 \subfigure[{\tiny  \(b<0\) and \(c<0\)} \label{cce}]
 {\includegraphics[width=.22\columnwidth,height=.22\columnwidth]{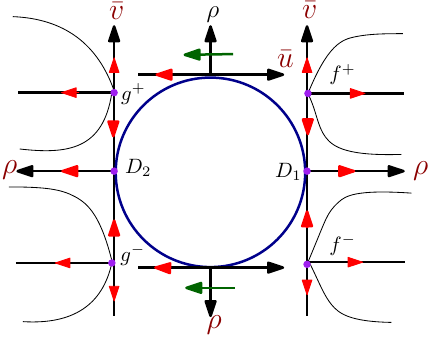}}
\quad
 \subfigure[{\tiny  \(b<0\) and \(c>0\)} \label{ccd}]
 {\includegraphics[width=.22\columnwidth,height=.22\columnwidth]{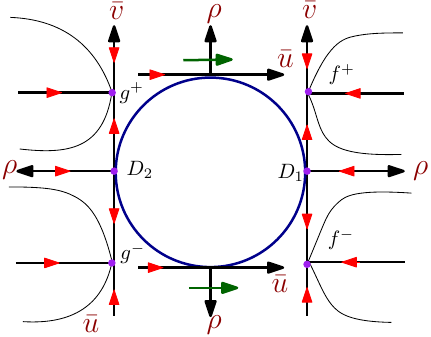}}
  \end{center}
  \caption{\label{cc} Qualitative properties of vector fields \(X_1^{\pm}\) and \(Y_1^{\pm}.\)  The axes corresponding to systems $X_1^{\pm}$ and $Y_1^{\pm}$ are $(\rho,\bar{v})$ and $(\rho,\bar{u}),$ respectively.  }
\end{figure}
\end{proof} 


 \section{Discussion and concluding remark}
\label{s: conclusion}

Although the FHN equations were originally developed as a simplified 
representation
of nerve impulses, they have also been extensively investigated 
in 
purely mathematical purposes
\cite{KS2001a, KS2001b, GLR2025, connections} 
 --  they consist of a very simple model of equations with  rich dynamics.

 Results obtained on  \cite[pp. 180]{fitzhughbook}  synthesised  the global bifurcation diagram for the FHN system \eqref{2d}. 
 It has been obtained by putting together, as in a ``huge puzzle'', all local bifurcation diagrams obtained. 
 The authors have concentrated their attention to  a particular parameter region relevant to physiology ($|c|>1+\sqrt{3}$). 

Trying to complete and understand the bifurcation diagram of \eqref{2d}, in this paper we have discussed the finite equilibria of \eqref{2d}, as well as their bifurcations, for 
 the  three scenarios identified in Subsection \ref{ss:structure}: 
\begin{description}
\item[Case A] \(a=0,\)
\item[Case B] \(b=0,\)
\item[Case C] \(a \neq 0$ and $ 0<b<1\). 
 \end{description}
 
The   richest scenario is \textbf{Case A}. In Theorem \ref{th:main1}  we have found   a Double-zero Bifurcation with a $\mathbb{Z}_2(\kappa)$--symmetry.   We have obtained precise expressions of the bifurcation curves passing through the bifurcation points $(b,c)=(1, \pm 1)$,  complementing the work started in \cite{GlobalFN}. 
We also give an analytical proof of the results stated in Section 2.4 of \cite{RS_proceedings}.

We have been able to explain rigorously the dynamics of Regions \textbf{8, 11, 17} and \textbf{18} of \cite{fitzhughbook} on the line defined by $a=0$ and $b>0$, as we proceed to explain:
\begin{itemize}
\item Point $\textbf{Q}$: Pitchfork bifurcation of $E_1$ (Theorem \ref{th:main1}, curve $b=T_P$); 
\item Region \textbf{8}: existence of three equilibria, one saddle $E_1$ and two sources $E_2, E_3$;
\item Point $\textbf{Q}_0$:  Hopf bifurcations of $E_2, E_3$ (Theorem \ref{th:main1}, curve $b=T_H^{2,3})$;
\item Region \textbf{11}: three periodic solutions (one stable, two unstable);
\item Point $\textbf{Q}_6$:  Double Homoclinic (Theorem \ref{th:main1}, curve $b= {DH}$); 
  \item Region \textbf{17}: two periodic solutions of different stabilities;
\item Point $\textbf{T} $:  Saddle-node bifurcation of non-hyperbolic solutions, making part of the Bautin bifurcation of $\textbf{Q}_{17}$ and $\textbf{Q}_{18}$ -- see Remark \ref{rem:SNL};
 \item Region \textbf{18}: three equilibria $E_1, E_2$ and $E_3$.
\end{itemize}
The  precise location the  above points/regions depend on $c $ but their relative position does not.  \\
 
The non-hyperbolic equilibria 
associated with
parameter values located at $\textbf{Q}$ (see \cite{fitzhughbook}) are degenerated saddle-node (cusp) equilibria, which are  attracting for $c < 1$ and repelling for $c > 1$. For $c = 1$, following the Remark   \ref{rem:SNL},  this point might correspond to a degenerated Bogdanov-Takens of order two  (with symmetry) -- a kind of codimension-three bifurcation. As pointed out in \cite{connections}, this would be the unique generic codimension-three local bifurcation exhibited by the FHN model \eqref{2d} and  its complete understanding is an open problem.  Using the same line of argument, in Cases \textbf{B} and \textbf{C}, we have analysed the dynamics of \eqref{2d}. 

In Section \ref{s:canards}, we have  studied the asymptotic case $c \to \pm \infty$, where canards are detected.
  A key point in this research
is the empirical difficulty 
to detect \emph{canards}
 because
 the smaller value of $\varepsilon$ is, the narrower the interval of $b$ (\textbf{Case A}) or $a$ (\textbf{Case B}) values for which \emph{canards} appear -- see Figure \ref{Axis}.  The dynamics that emerges in the system within these narrow parameter ranges is known as the \emph{{canard explosion}} 
and has been studied in \cite{GLR2025, Kuehn, KS2001a}.
 
 In Section \ref{sec4}, we have continued the analysis with the  compactification of the phase portraits of  \eqref{2d} on the Poincar\'e disc. This brings additional information of the trajectories which come from or tend to infinity.  We have provided   phase and bifurcation  diagrams in all the three cases. 
 
Throughout this article, we have connected the local bifurcation theory   (Theorems \ref{th:main1}, \ref{Th3.8} and Proposition \ref{Prop3.10}) with dynamics ``near infinity'' (Theorem \ref{S50_th}). Our contributions do not finish the whole discussion of the bifurcation analysis of \eqref{2d}; bifurcations might make part  of high codimension phenomena.

Theoretically, the maximum number of non-trivial periodic solutions of \eqref{2d} is three \cite[pp. 218]{fitzhughbook}. There are regions where the dynamics is completely determined and the phase portrait is complete and others in which the complete analysis of the  bifurcation diagram of  \cite[pp. 180]{fitzhughbook} is still ongoing.

\section{Appendix}

 All phase portrait  are illustrated with numerical simulations in Tables \ref{QA} and \ref{QB} using the software \emph{pplane9}.  
Initial conditions were rigorously chosen   to illustrate the characteristic behaviour in each region.
Furthermore, limit cycles are indicated by the red trajectories in the phase plot.
The phase portraits of Table \ref{QA} and \ref{QB} have been plotted for the parameter values, 
\begin{align*}
(c,b) =& (2.5,1.1),\; (2,1.2),\; (2,1.3),\; (2,1.34),\;
(1.3,1.3),\; (0.5,1.1),\; (-0.5,1.1),\; (-2,1.5),\\\nonumber
&(-2,1.34),\; (-2,1.33),\; (-2,1.2),\; (-2.5,1.1),\;
(-3,0.5),\; (-1,0.5),\; (-0.5,0.4)\;\\\nonumber
&(-0.1,0.9), \; (0.1,0.9),\; (0.4,0.4),\; (1.5,0.5),\; (3,0.5),\;
(2.5,-0.5),\; (2,-0.5),\;\\\nonumber
&(1,-1),\; (0.2,-1),\; (-0.2,-1),\; (-1,-1),\; (-1.5,-0.5),\; (-2.5,-0.5),
\end{align*}
and
\begin{align*}
(c,a) =& (1, 0.2),\; (-1,0.2),\; (5, 2),\; (-5, 2),\; (-1, 0.5),\; (1, 0.5),\; (-5, 0.5),\; (5, 0.5),
\end{align*}
respectively.
 
\newpage
\begin{center}
\resizebox{\textwidth}{!}{
\setlength{\tabcolsep}{12pt}
\begin{tabular}{c c c c c c c}
\parbox{0.12\textwidth}{\centering\rotatebox{90}{\includegraphics[width=6cm,height=4.5cm,keepaspectratio]{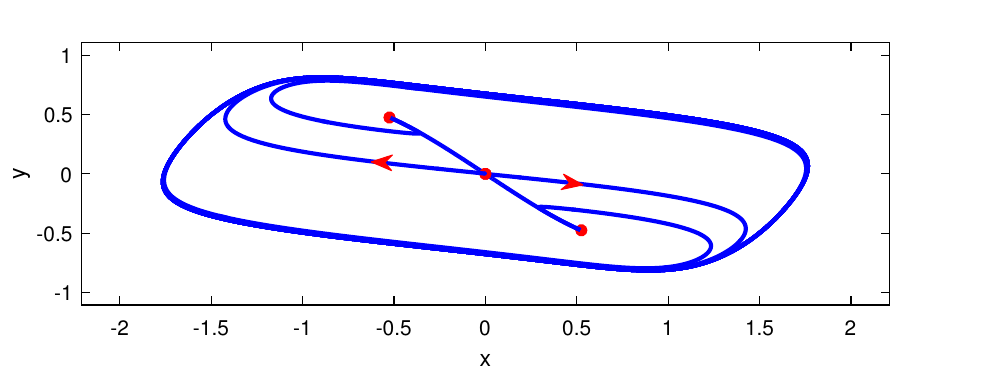}}\\ Region 1} &
\parbox{0.12\textwidth}{\centering\rotatebox{90}{\includegraphics[width=6cm,height=4.5cm,keepaspectratio]{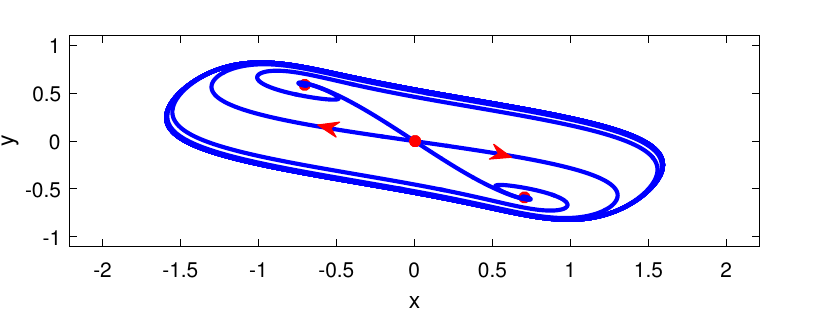}}\\ Region 2} &
\parbox{0.12\textwidth}{\centering\rotatebox{90}{\includegraphics[width=6cm,height=4.5cm,keepaspectratio]{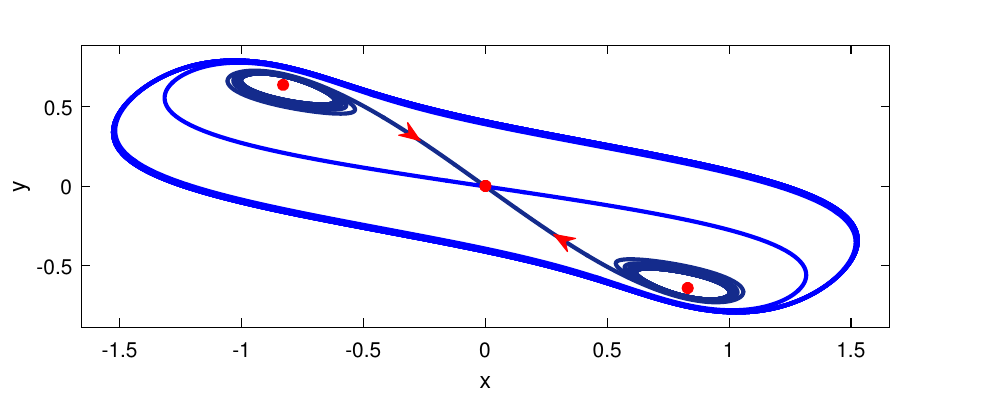}}\\ Region 3} &
\parbox{0.12\textwidth}{\centering\rotatebox{90}{\includegraphics[width=6cm,height=4.5cm,keepaspectratio]{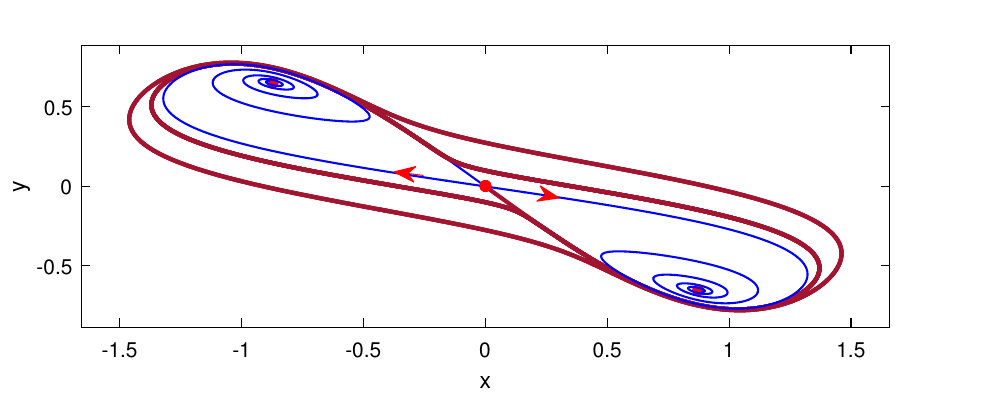}}\\ Region 4} &
\parbox{0.12\textwidth}{\centering\rotatebox{90}{\includegraphics[width=6cm,height=4.5cm,keepaspectratio]{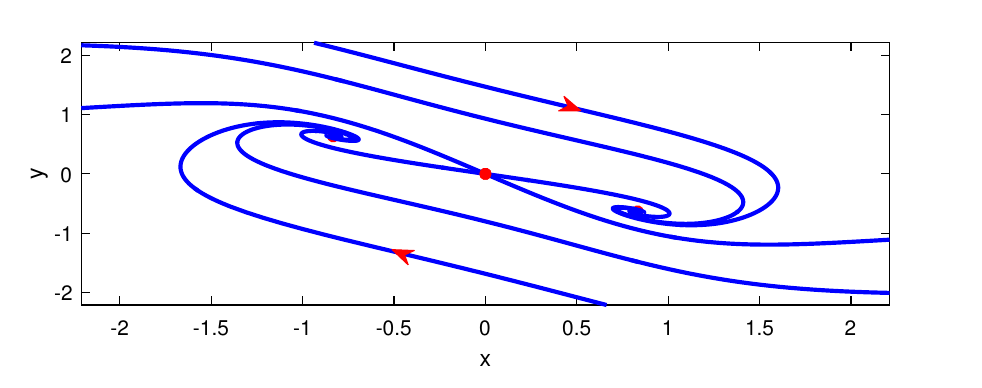}}\\ Region 5} &
\parbox{0.12\textwidth}{\centering\rotatebox{90}{\includegraphics[width=6cm,height=4.5cm,keepaspectratio]{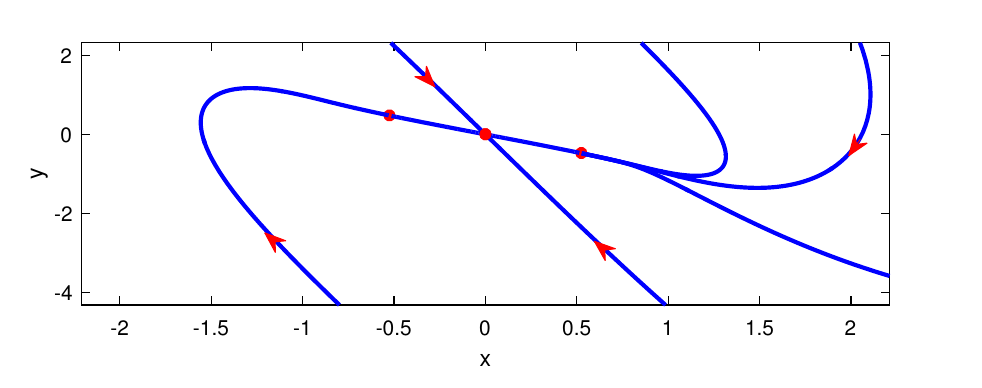}}\\ Region 6} &
\parbox{0.12\textwidth}{\centering\rotatebox{90}{\includegraphics[width=6cm,height=4.5cm,keepaspectratio]{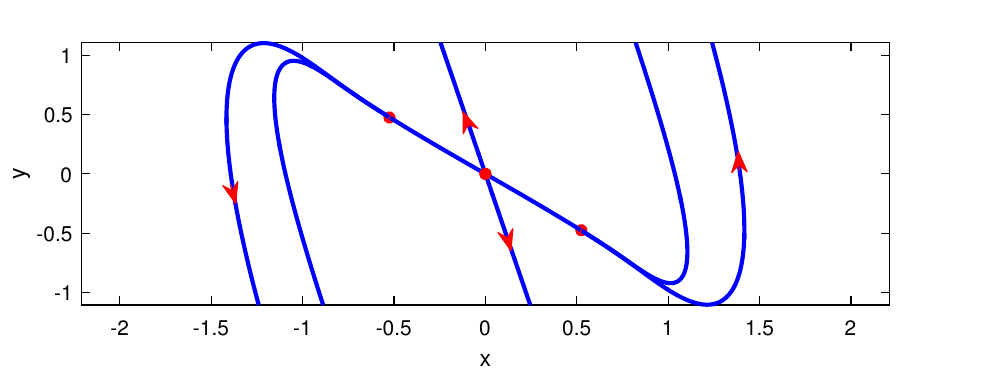}}\\ Region 7} \\
\parbox{0.12\textwidth}{\centering\rotatebox{90}{\includegraphics[width=6cm,height=4.5cm,keepaspectratio]{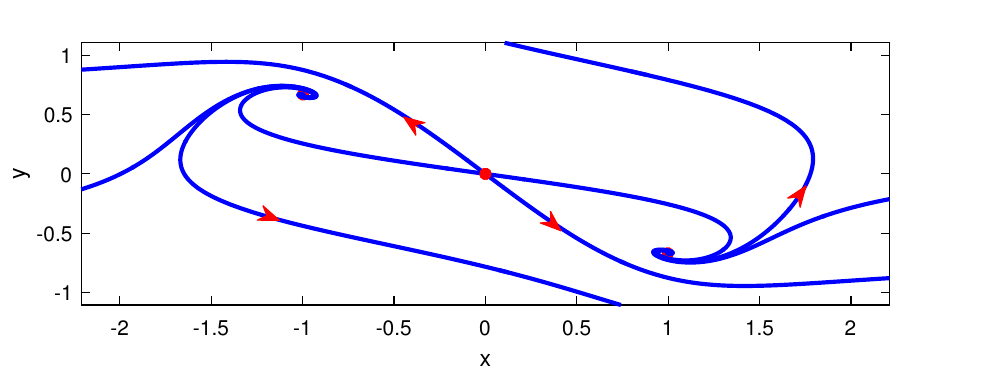}}\\ Region 8} &
\parbox{0.12\textwidth}{\centering\rotatebox{90}{\includegraphics[width=6cm,height=4.5cm,keepaspectratio]{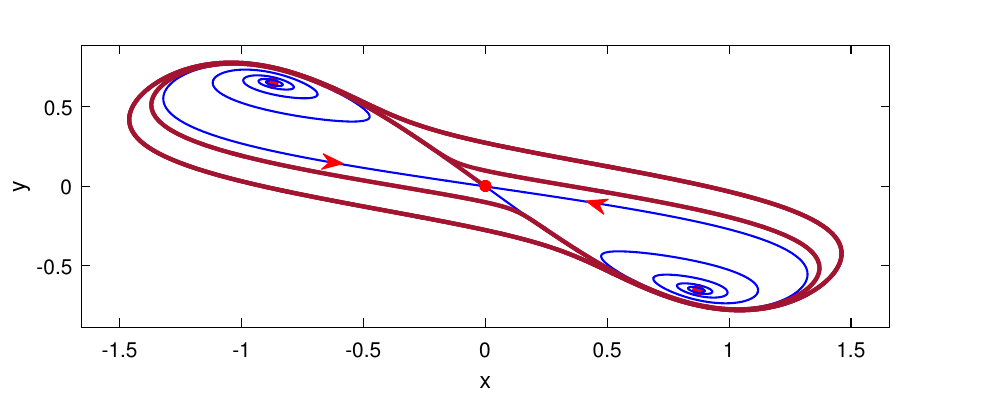}}\\ Region 9} &
\parbox{0.12\textwidth}{\centering\rotatebox{90}{\includegraphics[width=6cm,height=4.5cm,keepaspectratio]{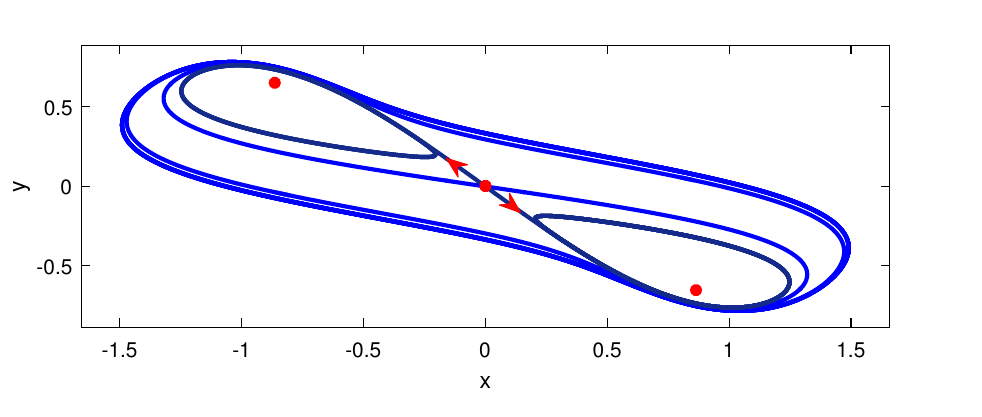}}\\ Region 10} &
\parbox{0.12\textwidth}{\centering\rotatebox{90}{\includegraphics[width=6.2cm,height=4.5cm,keepaspectratio]{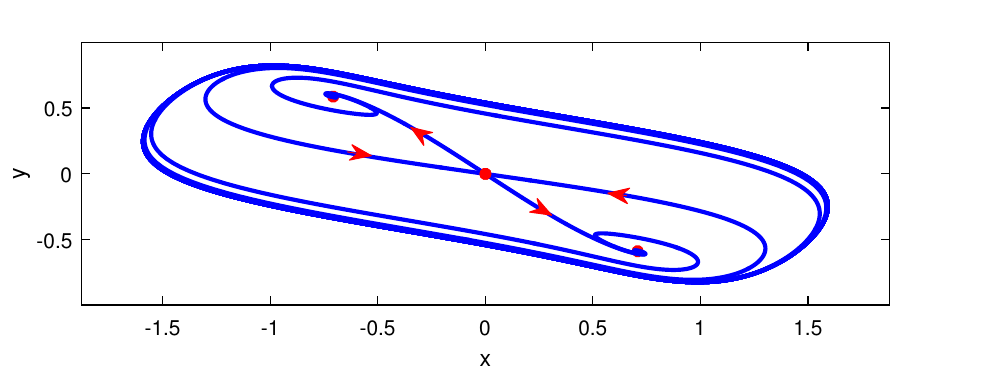}}\\ Region 11} &
\parbox{0.12\textwidth}{\centering\rotatebox{90}{\includegraphics[width=6.2cm,height=4.5cm,keepaspectratio]{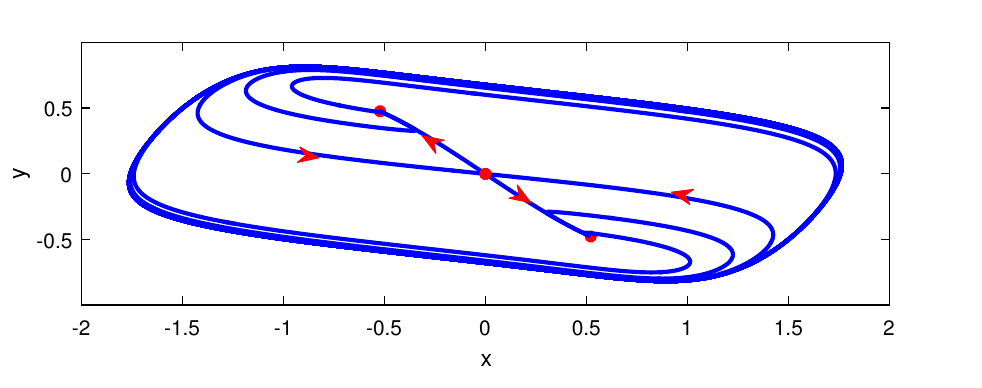}}\\ Region 12} &
\parbox{0.12\textwidth}{\centering\rotatebox{90}{\includegraphics[width=6cm,height=4.5cm,keepaspectratio]{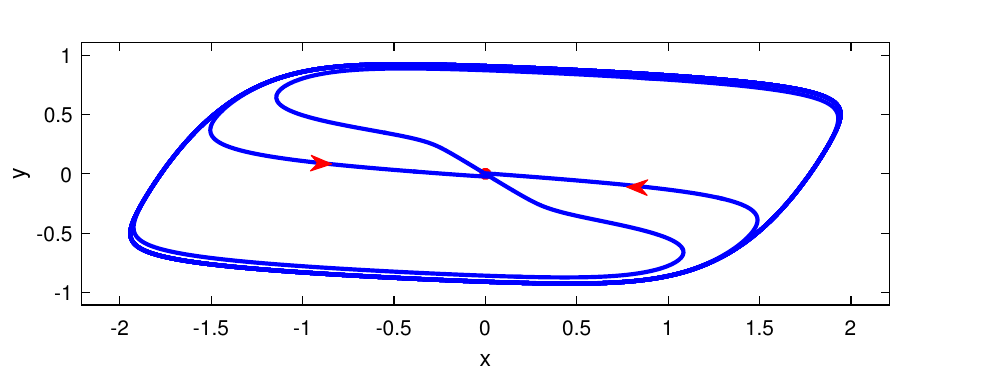}}\\ Region 13} &
\parbox{0.12\textwidth}{\centering\rotatebox{90}{\includegraphics[width=6cm,height=4.5cm,keepaspectratio]{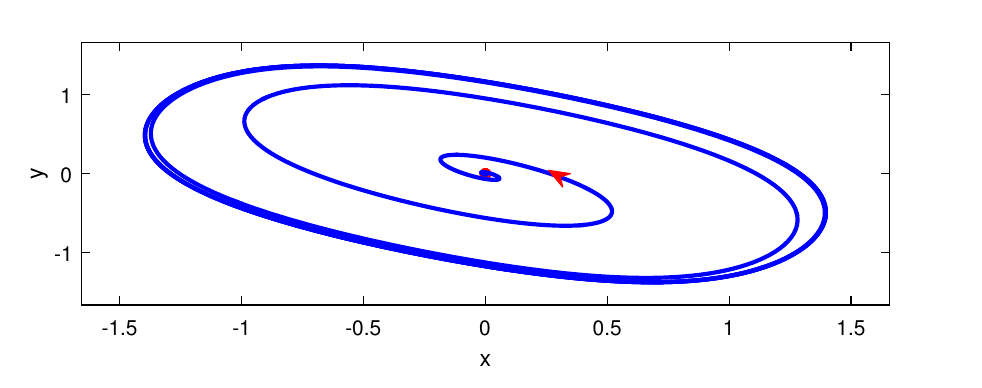}}\\ Region 14} \\ 
\parbox{0.12\textwidth}{\centering\rotatebox{90}{\includegraphics[width=6cm,height=4.5cm,keepaspectratio]{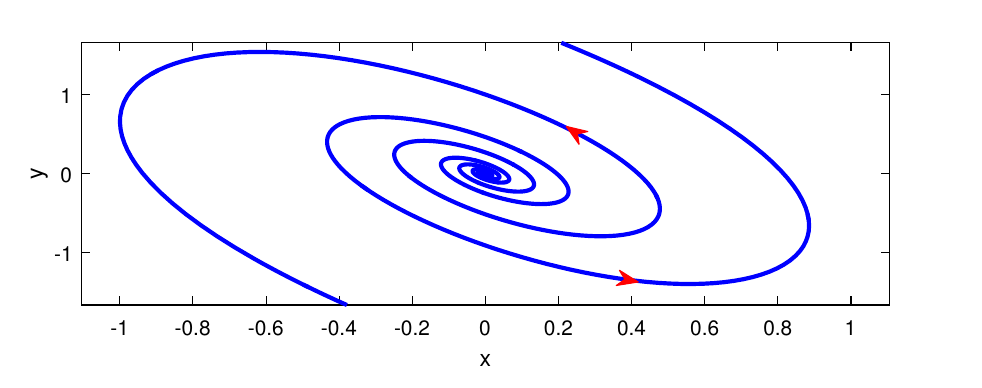}}\\ Region 15} &
\parbox{0.12\textwidth}{\centering\rotatebox{90}{\includegraphics[width=6cm,height=4.5cm,keepaspectratio]{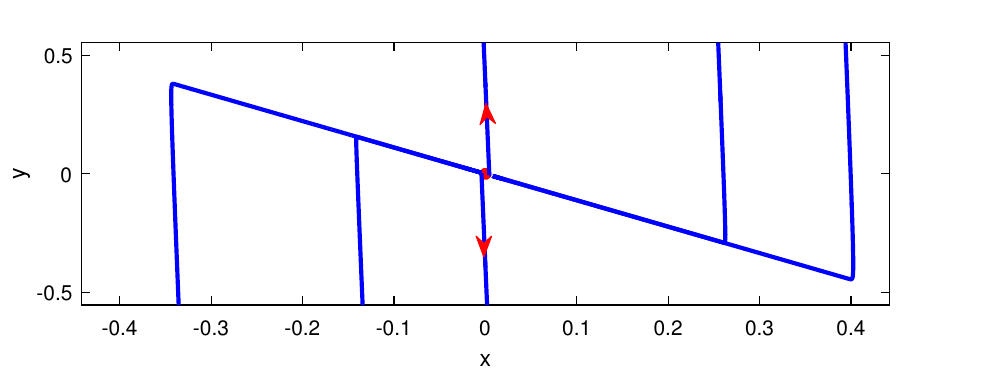}}\\ Region 16} &
\parbox{0.12\textwidth}{\centering\rotatebox{90}{\includegraphics[width=6cm,height=4.5cm,keepaspectratio]{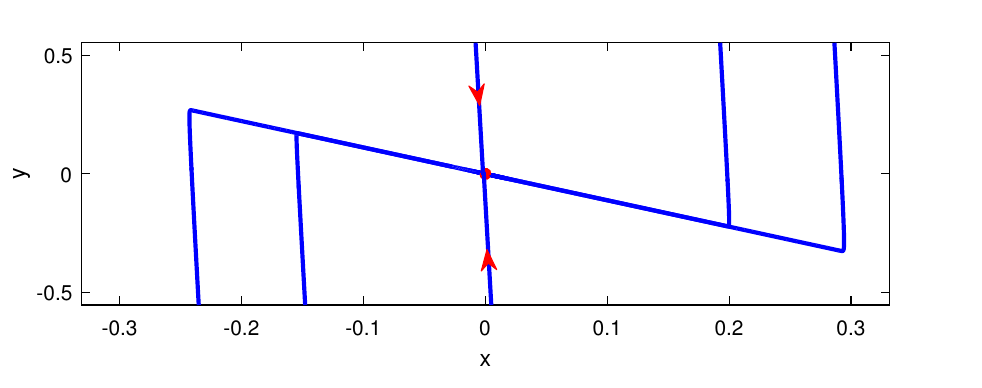}}\\ Region 17} &
\parbox{0.12\textwidth}{\centering\rotatebox{90}{\includegraphics[width=6cm,height=4.5cm,keepaspectratio]{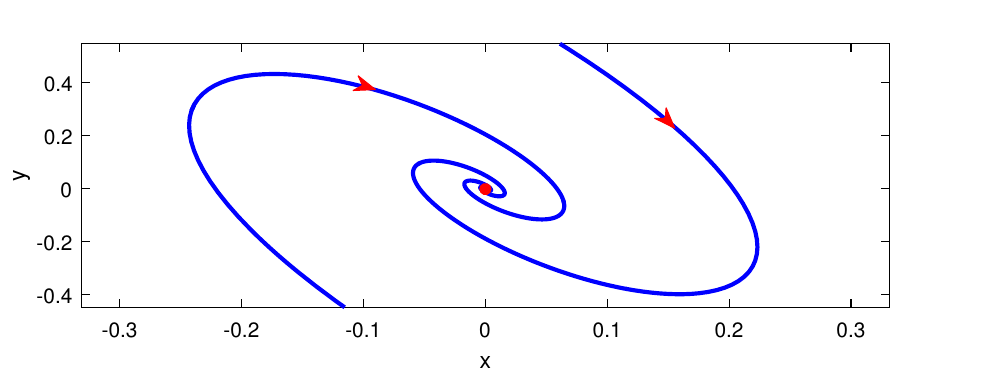}}\\ Region 18} &
\parbox{0.12\textwidth}{\centering\rotatebox{90}{\includegraphics[width=6cm,height=4.5cm,keepaspectratio]{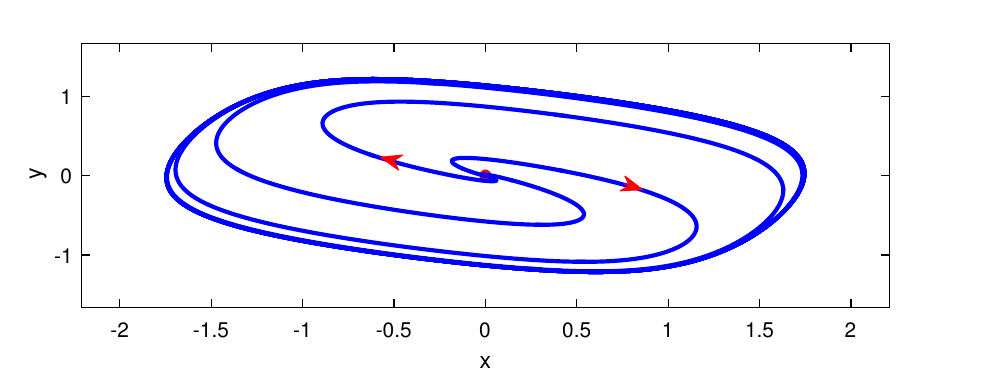}}\\ Region 19} & 
\parbox{0.12\textwidth}{\centering\rotatebox{90}
{\includegraphics[width=6cm,height=4.5cm,keepaspectratio]{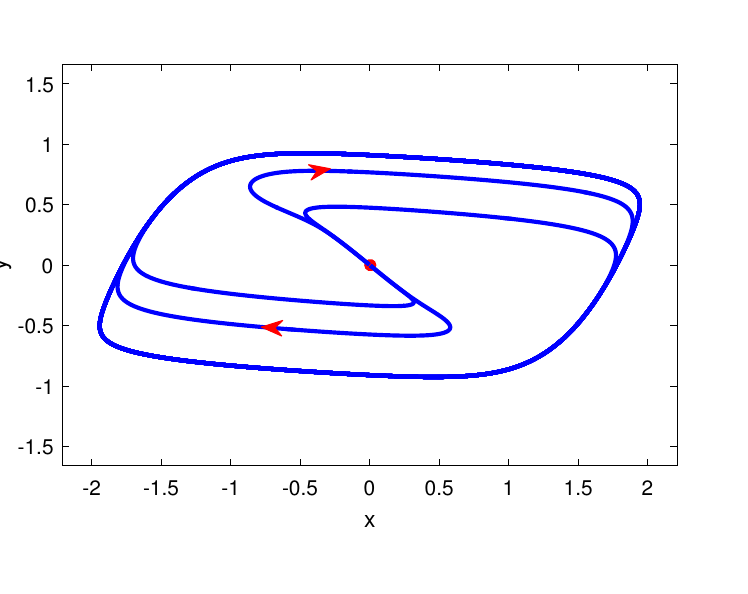}}\\ Region 20} &
\,\,\,
\parbox{0.12\textwidth}{\centering\rotatebox{90}{\includegraphics[width=6cm,height=4.5cm,keepaspectratio]{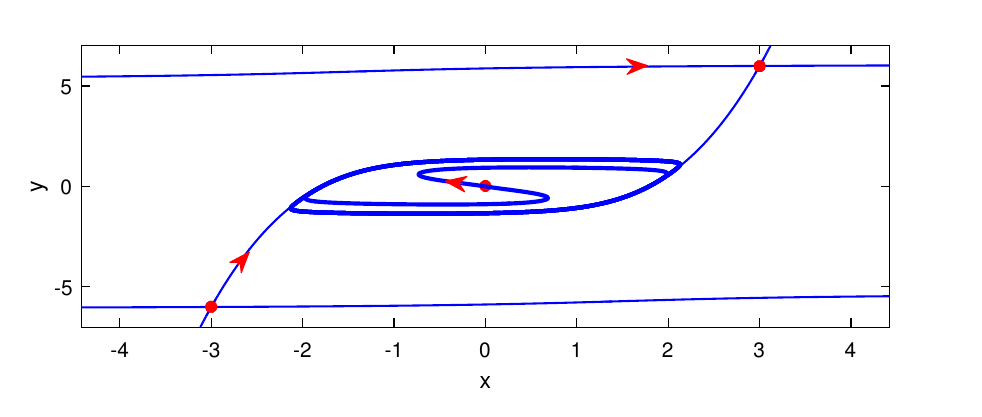}}\\ Region 21} \\
\parbox{0.12\textwidth}{\centering\rotatebox{90}{\includegraphics[width=6cm,height=4.5cm,keepaspectratio]{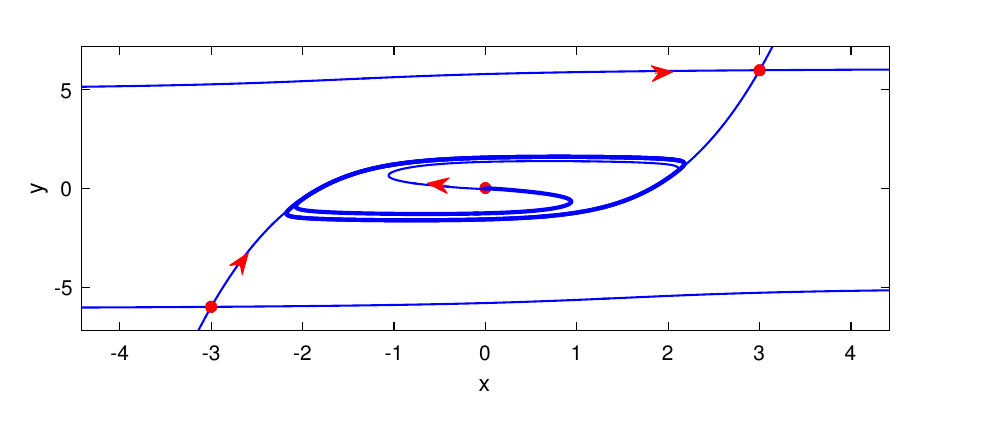}}\\ Region 22} &
\parbox{0.12\textwidth}{\centering\rotatebox{90}{\includegraphics[width=6cm,height=4.5cm,keepaspectratio]{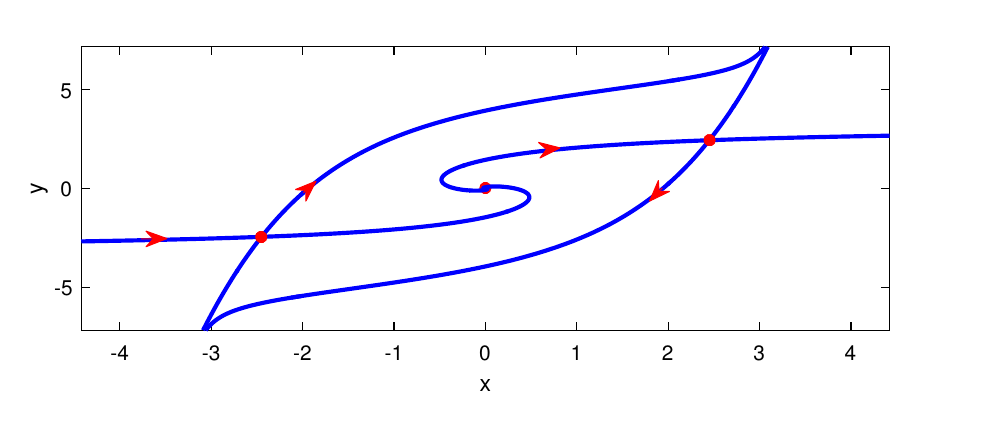}}\\ Region 23} &
\parbox{0.12\textwidth}{\centering\rotatebox{90}{\includegraphics[width=6cm,height=4.5cm,keepaspectratio]{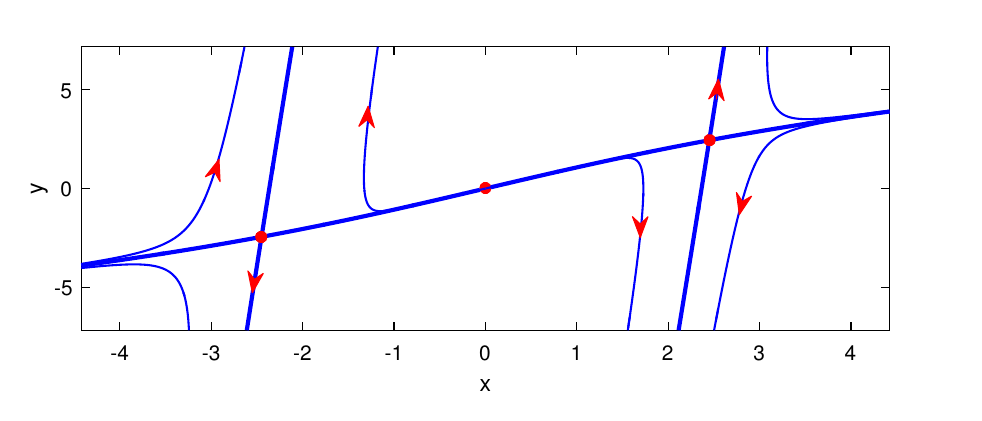}}\\ Region 24} &
\parbox{0.12\textwidth}{\centering\rotatebox{90}{\includegraphics[width=6cm,height=4.5cm,keepaspectratio]{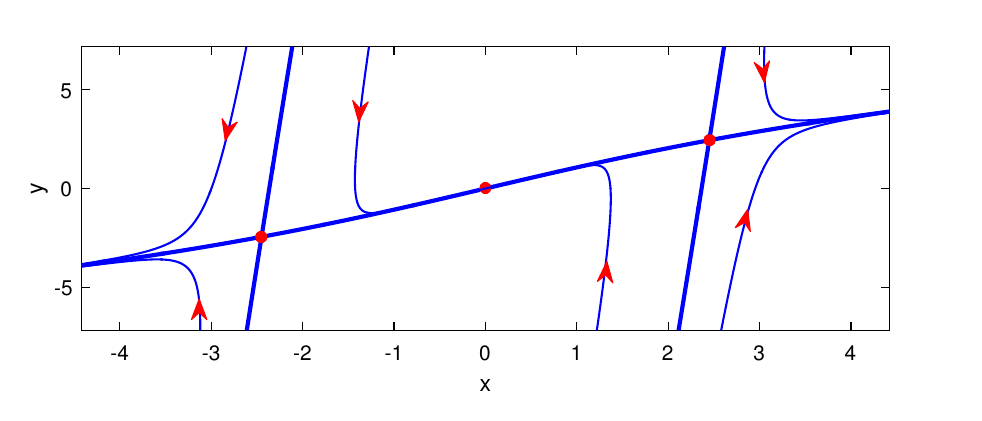}}\\ Region 25} &
\parbox{0.12\textwidth}{\centering\rotatebox{90}{\includegraphics[width=6cm,height=4.5cm,keepaspectratio]{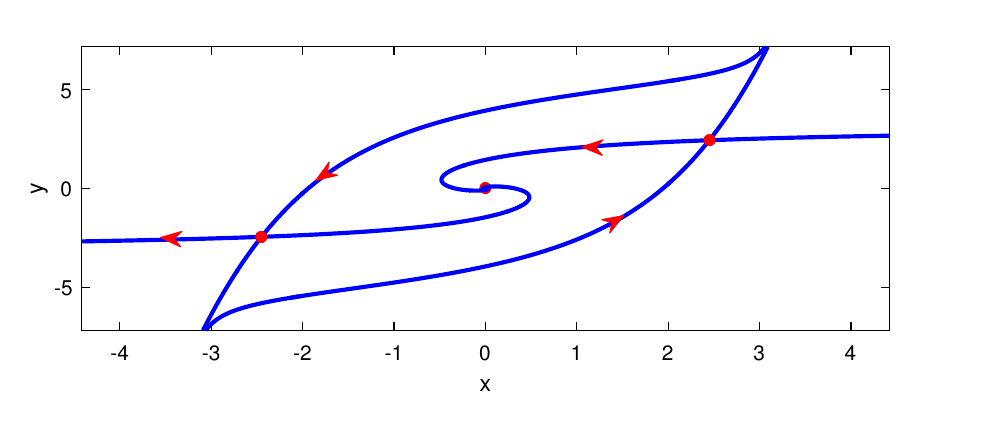}}\\ Region 26} &
\parbox{0.12\textwidth}{\centering\rotatebox{90}{\includegraphics[width=6cm,height=4.5cm,keepaspectratio]{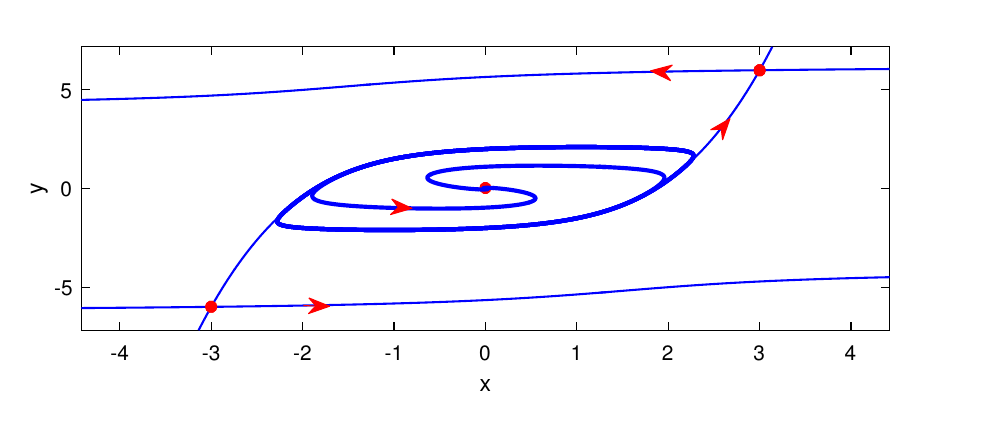}}\\ Region 27} &
\parbox{0.12\textwidth}{\centering\rotatebox{90}{\includegraphics[width=6cm,height=4.5cm,keepaspectratio]{e24new}}\\ Region 28} \\
\end{tabular}
}
\captionof{table}{\label{QA} The local phase portraits corresponding to \textbf{Case A}, which is associated with system \eqref{2d} when \(a=0,\) for Regions \textbf{1--28} of Figure \ref{CASEA}.}

\end{center}

\begin{center}
\begin{table}[H]
\setlength{\tabcolsep}{10pt}
\renewcommand{\arraystretch}{1.5}
\begin{tabular}{cccc}
{\includegraphics[width=.23\columnwidth,height=.20\columnwidth]{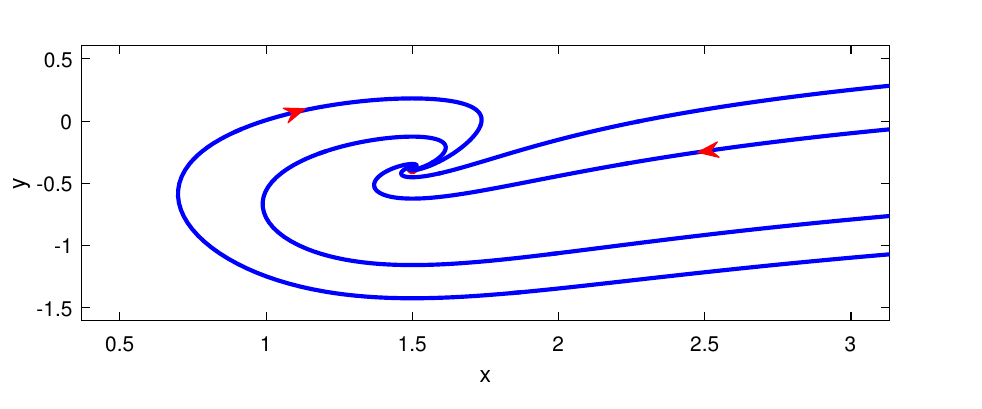}}&
{\includegraphics[width=.23\columnwidth,height=.20\columnwidth]{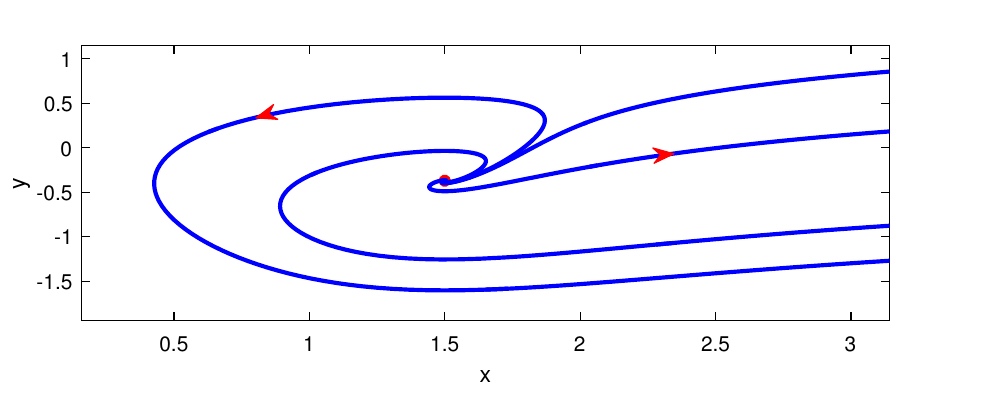}}&
{\includegraphics[width=.24\columnwidth,height=.20\columnwidth]{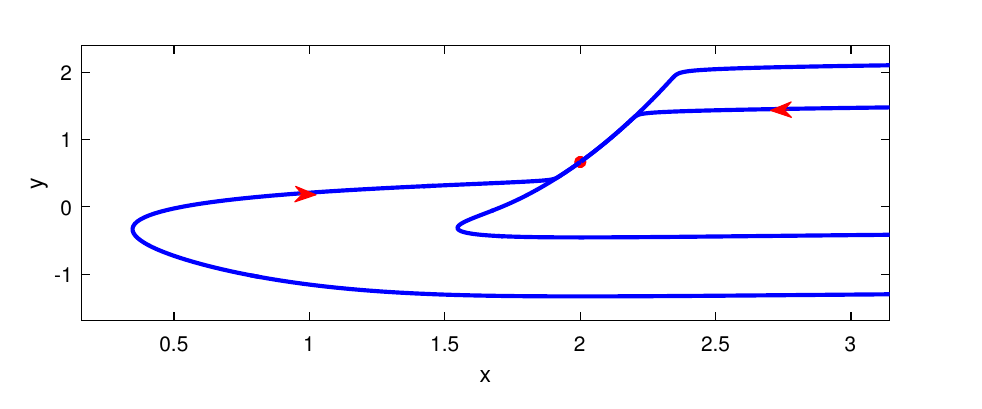}}&
{\includegraphics[width=.24\columnwidth,height=.20\columnwidth]{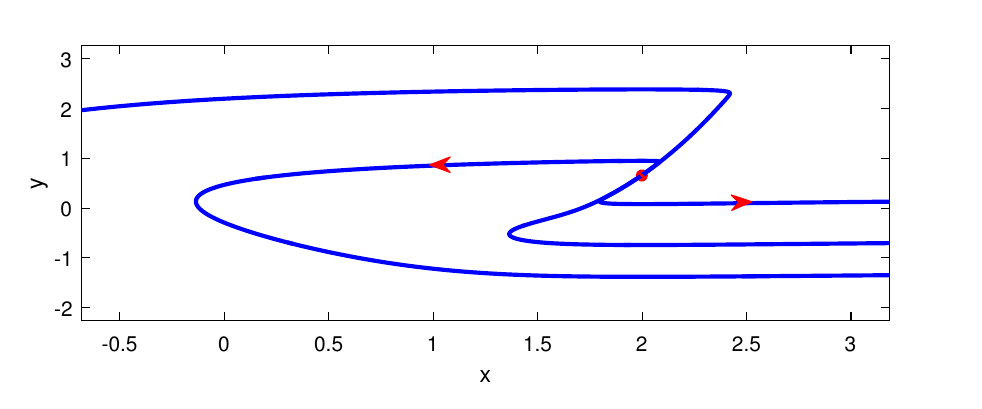}}\\
Region 1 & Region 2 & Region 3 & Region 4 \\
{\includegraphics[width=.23\columnwidth,height=.20\columnwidth]{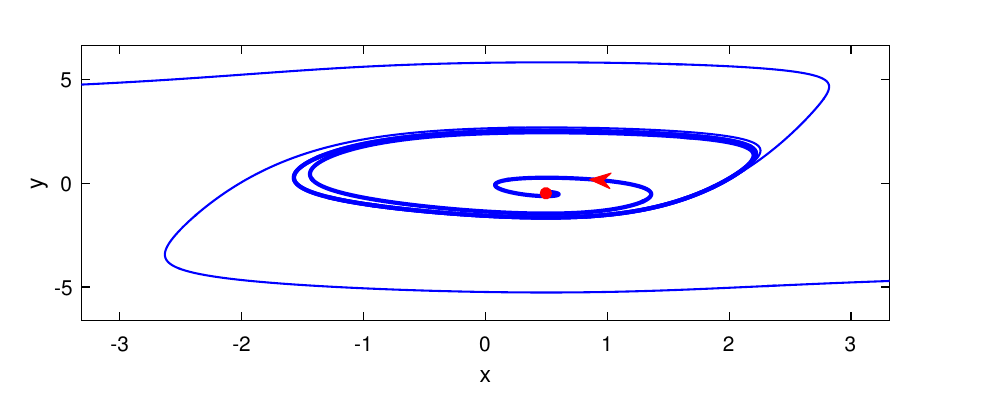}}&
{\includegraphics[width=.23\columnwidth,height=.20\columnwidth]{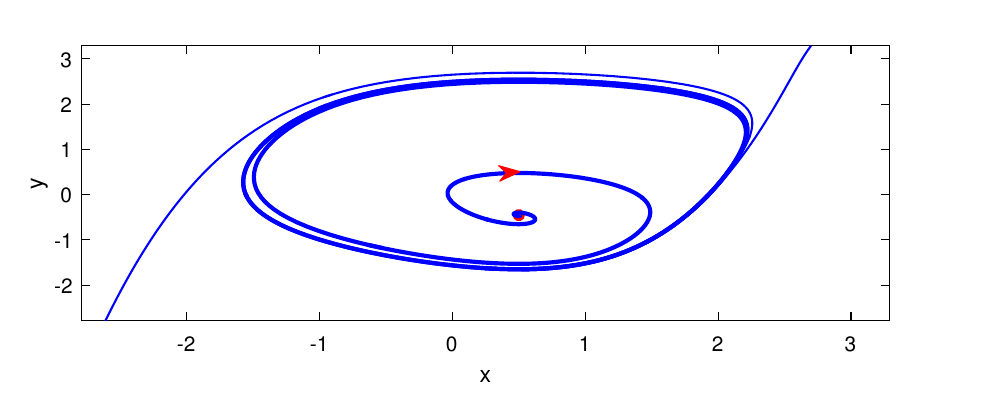}}&
{\includegraphics[width=.24\columnwidth,height=.20\columnwidth]{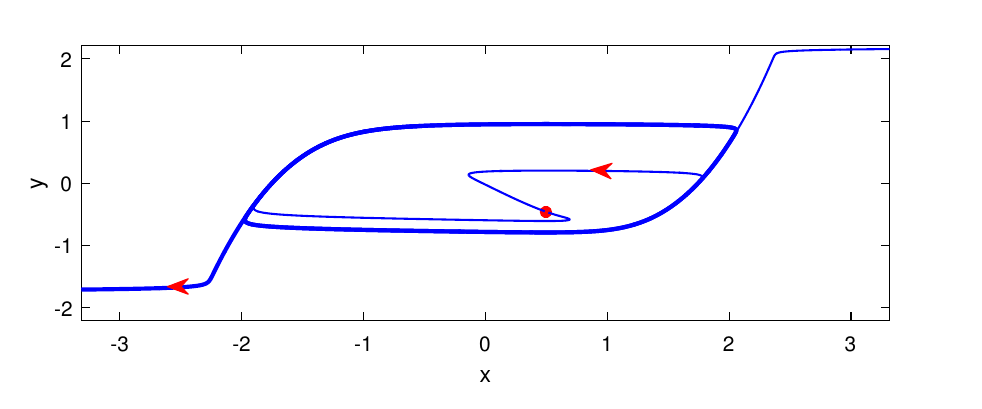}}&
{\includegraphics[width=.24\columnwidth,height=.20\columnwidth]{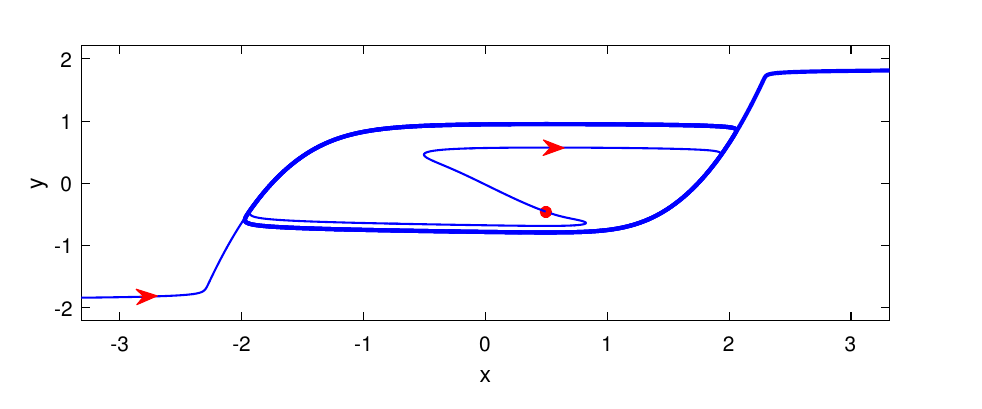}}\\
Region 5 & Region 6 & Region 7 & Region 8 \\
\end{tabular}
\caption{\label{QB} The local phase portraits corresponding to \textbf{Case B}, which is associated with system \eqref{2d} when \(b=0,\) for Regions \textbf{1--8} of Figure \ref{CASEC}. }
\end{table}
\end{center}

\section*{Aknowledgments}
The authors are grateful to Professor Jaume Llbre for the formal proof of Theorem \ref{thm11}.


\begin{thebibliography}{99}

\bibitem{Alonso}
Sergio Alonso, Markus B\"ar and Blas Echebarria, 
{\em Nonlinear physics of electrical wave propagation in the heart: a review,}
Reports on Progress in Physics, {\bf 79}(9) (2016), 096601.

\bibitem{AFJ} M.J. Alvarez, A. Ferragut and X. Jarque,
A survey on the blow-up technique,  Internat. J. Bifur. Chaos Appl. Sci. Engrg. {\bf 21}  (2011), 3103--3118.


\bibitem{CPRG2024} D. Cebri\'an-Lacasa, P. Parra-Rivas, D. Ruiz-Reyn\'es, L. Gelens, 
{\sl Six decades of the FitzHugh--Nagumo model:
A guide through its spatio-temporal dynamics and influence across disciplines} 
Physics Reports, 1096, 1--39, 2024.


\bibitem{cherkas}
   L.A. Cherkas and L.I. Zhilevich,
    Some tests for the absence or uniqueness of limit cycles,
    \textit{Differentsial'nye Uravneniya},
    6 , 1170--1178, 1970.

\bibitem{Abid} N. Ben Abid, M. Bendahmane, M. Mahjoub, 
{\em Stability of the Ionic Parameters of a Nonlocal FitzHugh-Nagumo Model of Cardiac Electrophysiology,} 
Acta Applicandae Mathematicae, {\bf 193}(1), 2024

\bibitem{llibrebook}
F. Dumortier, J.  Llibre, J.C. Art\'es, 
Qualitative theory of planar differential systems,
Springer, New York, 2006.

\bibitem{Fenichel}
N. Fenichel, (1971). {\em Persistence and smoothness of invariant manifolds for flows,}
Indiana University Mathematics Journal, {\bf 21} (3), 193--226.

\bibitem{Fitzhugh}
R. Fitzhugh, 
{\em Impulses and physiological states in theoretical models and propagation in nerve membrane,}
Biophys. J., {\bf 1} (445), 1961.

\bibitem{GlobalFN}
Adelina Georgescu, Carmen Roc\c soreanu, and Nicolaie Giurgi\c teanu,
{\em Global Bifurcations in FitzHugh-Nagumo Model}, Trends in Mathematics:
Bifurcations, Symmetry and Patterns,  197-202, 2003.

 \bibitem{GS85}
M. Golubitsky, D.G. Schaeffer, Singularities and groups in bifurcation theory, Vol. 1,
Springer-Verlag , 1985.

\bibitem{GLR2025} 
B.F.F. Gon\c{c}alves, I.S. Labouriau, A.A.P. Rodrigues, \emph{Bifurcations and canards in the FitzHugh-Nagumo system: a tutorial of fast-slow dynamics}, International Journal of Bifurcation and Chaos, Vol. 35 (8), 2350017, 2025.


\bibitem{HoHu}
A. L. Hodgkin,  A. F. Huxley, 
{\em A quantitative description of membrane current and its application to conduction and excitation
in nerve,} J. Physiol., {\bf 117},  500--544, 1952.

\bibitem{KS2001a}
M. Krupa,  P. Szmolyan,
\emph{Extending geometric singular perturbation theory to nonhyperbolic
points--fold and canard points in two dimensions}, SIAM Journal on Mathematical Analysis 33, 2001.
 
 \bibitem{KS2001b}
M. Krupa,  P. Szmolyan,
\emph{Relaxation oscillation and canard explosion}, Journal of Differential Equations 174, 312--368, 2001.

\bibitem{Kuehn} C. Kuehn, {\sl Multiple time scale dynamics}, Springer-Verlag 2015.

\bibitem{Kuznetsov}
Yuri A. Kuznetsov,  Elements of Applied Bifurcation Theory. 3rd ed., Springer, 2004.


\bibitem{Linard}
A. Li\'enard, {\em  Etude des oscillations entretenues,} Revue G\'en\'erale de lElectricit\'e {\bf 23},  901--912 and 946--954, 1928.


%



\bibitem{Meiss} 
James. D. Meiss, Differential Dynamical Systems, SIAM, Philadelphia, PA, 2007.

\bibitem{CR2023} 
J. P. Maur\'icio de Carvalho, A. Rodrigues, \emph{SIR model with vaccination: bifurcation analysis},  Qualitative theory of dynamical systems 22(3), 105, 2023.


\bibitem{Nagumo}   
J. Nagumo, S. Arimoto, S. Yoshizawa,
{\em An active pulse transmission line simulating nerve axon,}
Proceedings of the IRE, {\bf 50}, 2061--2070. 1963.



\bibitem{RS_proceedings} C. Rocsoreanu, M. Sterpu, \emph{Local Bifurcation for the Fitzhugh-Nagumo System},  IFIP International Information Security Conference. Boston, MA: Springer US, 2002.

\bibitem{fitzhughbook}
C. Rocsoreanu, A. Georgescu, N. Giurgiteanu, FitzHugh-Nagumo Model: Bifurcation and Dynamics, 2nd ed., Kluwer Academic Publishers, Dordrecht, 2000.

\bibitem{2D}
C. Rocsoreanu, N. Giurgiteanu, A. Georgescu,
{\em Degenerated Hopf bifurcation in the Fitzhugh-Nagumo system. 2. Bautin bifurcation,}
Revue D\'analyse Num\'erique et de Th\'eorie de L\'approximation,
{\bf 29}(1),  97--109, 2000.

\bibitem{connections}
C. Rocsoreanu, N. Giurgiteanu, A. Georgescu,
{\em Connections between saddles for the FitzHugh-Nagumo system.} Int. J. Bif. and Chaos, 
{\bf 11}(2),   533--540, 2001. 


\bibitem{Tsai}
Tony Yu-Chen Tsai, Yoon Sup Choi, Wenzhe Ma, Joseph R Pomerening, Chao Tang, James E Ferrell Jr,
{\em Robust, tunable biological oscillations from interlinked positive and negative feedback loops,}
Science, {\bf 321}(5885),  126-9, 2008

\bibitem{vander}
B. van der Pol, 
{\em A theory of the amplitude of free and forced triode vibrations,} Radio Review, {\bf 1}, 701--710, and 754--762, 1920.

\bibitem{Verduzco}
G. Fernando Verduzco, 
{\em The first Lyapunov coefficient for a class of systems,} IFAC Proceedings Volumes, {\bf 38}(1) ,  1205--1209, 2005.

\bibitem{wiggins}
Stephen Wiggins, Introduction to Applied Nonlinear Dynamical Systems and Chaos. 2nd ed., Springer, 2003.

\end{thebibliography}
\end{document}